\documentclass[12pt,reqno]{article}
mm \textheight=8.9in

\usepackage{amssymb}
\usepackage{amsmath}
\usepackage{amsthm}

\usepackage{color}
\usepackage{pb-diagram}
\newcommand{\tw}[3]{{$#1$}${\,\scriptscriptstyle {#2}}\atop\raise9pt\hbox{$\scriptstyle\tp$} ${$#3$}}
\newcommand{\st}[1]{\mbox{${\,\scriptscriptstyle {#1}}\atop\raise5.5pt\hbox{$*$}$}}
\newcommand{\btr}{\raise1.2pt\hbox{$\scriptstyle\blacktriangleright$}\hspace{2pt}}

\newcommand{\id}{\mathrm{id}}
\newcommand{\im}{\mathrm{im}\:}

\newcommand{\A}{\mathcal{A}}

\newcommand{\Pc}{\mathcal{P}}

\newcommand{\D}{\mathfrak{D}}

\newcommand{\Sg}{\mathfrak{S}}
\newcommand{\Jg}{\mathfrak{J}}
\newcommand{\Tg}{\mathfrak{T}}

\newcommand{\Ru}{\mathcal{R}}

\newcommand{\E}{\mathcal{E}}

\newcommand{\Kc}{\mathcal{K}}
\newcommand{\Q}{\mathcal{Q}}
\renewcommand{\O}{\mathcal{O}}
\newcommand{\Nc}{\mathcal{N}}

\newcommand{\C}{\mathbb{C}}
\newcommand{\Z}{\mathbb{Z}}

\newcommand{\zz}{z}
\newcommand{\tp}{\otimes}
\newcommand{\vt}{\vartheta}

\newcommand{\U}{\mathcal{U}}

\newcommand{\ve}{\varepsilon}
\newcommand{\gm}{\gamma}

\newcommand{\la}{\lambda}
\newcommand{\tr}{\triangleright}

\newcommand{\ch}{\mathrm{ch}}

\newcommand{\End}{\mathrm{End}}

\newcommand{\Hom}{\mathrm{Hom}}

\newcommand{\Ind}{\mathrm{Ind}}
\newcommand{\rk}{\mathrm{rk}}
\newcommand{\Tr}{\mathrm{Tr}}

\newcommand{\Rm}{\mathrm{R}}
\newcommand{\Wm}{\mathrm{W}}

\newcommand{\btl}{\mbox{\raise1.1pt\hbox{$\scriptstyle\blacktriangleleft$}}}
\newcommand{\ad}{\mathrm{ad}}
\newcommand{\Ad}{\mathrm{Ad}}

\newcommand{\g}{\mathfrak{g}}
\renewcommand{\b}{\mathfrak{b}}

\newcommand{\h}{\mathfrak{h}}
\newcommand{\mub}{\boldsymbol{\mu}}

\newcommand{\n}{\mathfrak{n}}
\newcommand{\nb}{\boldsymbol{n}}
\newcommand{\m}{\mathfrak{m}}

\newcommand{\eps}{\epsilon}

\newcommand{\nn}{\nonumber}

\newcommand{\p}{\mathfrak{p}}
\renewcommand{\l}{\mathfrak{l}}
\renewcommand{\c}{\mathfrak{c}}

\newcommand{\al}{\alpha}

\newcommand{\bt}{\beta}
\newcommand{\be}{\begin{eqnarray}}
\newcommand{\ee}{\end{eqnarray}}

\newtheorem{thm}{Theorem}[section]
\newtheorem{propn}[thm]{Proposition}
\newtheorem{lemma}[thm]{Lemma}
\newtheorem{corollary}[thm]{Corollary}

\theoremstyle{definition}
\newtheorem{remark}[thm]{Remark}
\newtheorem{definition}[thm]{Definition}

\newcommand{\select}[1]{\textcolor{red}{\bf\em #1}}

\begin{document}
\title{Quantum conjugacy classes of simple matrix groups\footnote{
This research is partially supported
by
the Emmy Noether Research Institute for Mathematics,
the Minerva Foundation of Germany,  the Excellency Center "Group
Theoretic Methods in the study of Algebraic Varieties"  of the Israel
Science foundation, and by the RFBR grant no. 03-01-00593.
}}
\author{A. Mudrov
\\
\small
\sl Dedicated to the memory of Joseph Donin
}
\date{}
\maketitle
\begin{center}
{Emmy Noether Mathematics Institute,  52900 Ramat Gan,
Israel,\\
Max-Planck Institut f$\ddot{\rm u}$r Mathematik, Vivatsgasse 7, D-53111 Bonn, Germany.\\
e-mail: mudrova@macs.biu.ac.il, mudrov@mpim-bonn.mpg.de}
\end{center}
\begin{abstract}
Let $G$ be a  simple complex classical group and $\g$ its Lie algebra.
Let $\U_\hbar(\g)$ be the Drinfeld-Jimbo quantization
of the universal enveloping algebra $\U(\g)$.
We construct an explicit $\U_\hbar(\g)$-equivariant quantization
of conjugacy classes of $G$ with Levi subgroups as the stabilizers.
\end{abstract}
\maketitle
\tableofcontents

\section {Introduction}
Deformation quantization of a Poisson structure on a smooth manifold is a classic problem of mathematical physics.
Especially interesting is a quantization that is equivariant with respect to an action
of a group and,  even more generally, a quantum group. Recently, a significant progress
in this field was triggered by a discovered connection between equivariant quantization and the theory of dynamical
Yang-Baxter equation, \cite{DM1}. Namely, a star product quantization of semisimple coadjoint
orbits and conjugacy classes (with a Levi subgroup as the stabilizer) of simple Lie
groups was constructed in \cite{EE,EEM} in terms of
the universal dynamical twist.

On the other hand, semisimple coadjoint orbits and conjugacy classes
are affine algebraic varieties. Therefore their quantization may be sought for
in terms of generators and relations, as a deformed ring of polynomial functions. This is an alternative approach
as compared to star product and it has certain advantages because the solution is formulated
by "finite data".

In the present paper, we construct an equivariant deformation quantization of semisimple
conjugacy classes of simple classical algebraic
groups with Levi subgroups as the stabilizers.
More precisely, let $G$ be a complex simple algebraic group from the series $A$, $B$, $C$, and $D$.
Fix the standard, or Drinfeld-Jimbo, factorizable quasitriangular Lie bialgebra structure on $\g =\mathrm{Lie}\;G$.
Equipped with the corresponding Drinfeld-Sklyanin (DS) bracket, $G$ becomes a  Poisson Lie group. Consider
the adjoint action of $G$ on itself.
The Semenov-Tyan-Shansky (STS) Poisson bracket makes  $G$ a Poisson Lie manifold
over $G$. The symplectic leaves of this Poisson structure
are exactly the conjugacy classes.
Let  $\U_\hbar(\g)$ be the Drinfeld-Jimbo quantum group.
We construct a $\U_\hbar(\g)$-equivariant quantization of the ring of polynomial functions
along the STS bracket
on almost all conjugacy classes with Levi subgroups as the stabilizers. Almost all means
all for special linear and symplectic groups. For the orthogonal groups our construction
covers the
classes of matrices with eigenvalues $\la$ subject to the condition $\la^2=1\Rightarrow \la=1$. Those are exactly
the classes that are isomorphic to semisimple  $G$-orbits in $\g^*$ (via the Cayley transformation).

The quantization is given explicitly, in terms of deformed ideals of classes
in a quantized ring of polynomial functions on the group (or in the so called
reflection equation algebra).
Simultaneously, the quantized classes  are realized as subalgebras of operators
on generalized Verma modules.

The setup of the paper is as follows.

In Section \ref{secDJQG} we recall definition of the Drinfeld-Jimbo quantum group and its some
important subalgebras.

In Section \ref{secGVM} we study tensor products of finite dimensional and generalized Verma modules over
$\U_\hbar(\g)$.

In Section \ref{secREM} we study properties of a fundamental object of our theory, the universal
reflection equation matrix $\Q=\Ru_{21}\Ru$ (expressed through the universal R-matrix of $\U_\hbar(\g)$).
In particular, we determine its spectrum on
tensor products of finite dimensional and generalized Verma modules and compute
q-traces of $\Q^\ell$.

In Section \ref{secQAAV} we develop a method of quantization of affine homogeneous varieties.

In Section \ref{secQSAG} we recall some results of \cite{M} on  quantization of simple algebraic
groups and construct an embedding of the quantized affine coordinate ring $\C_\hbar[G]$ in $\U_\hbar(\g)$.

In Section \ref{subsecMC} we recall results of \cite{M} concerning the center of $\C_\hbar[G]$.

In Section \ref{secQCC} we give the quantization of the conjugacy classes.

Appendix contains auxiliary an information about invariants of $\U_\hbar(\g)$.

\vspace{0.5cm}
\noindent
{\bf \large Acknowledgements.}
The author is grateful to the Max-Planck Institute for Mathematics in Bonn for hospitality and
the excellent research conditions. He thanks J. Bernstein for his interest to the present work, helpful
discussions and valuable remarks.

\section{Drinfeld-Jimbo quantum group}
\label{secDJQG}
In the present paper, we work over the ring $\C[[\hbar]]$
of formal power series in $\hbar$.
Given a $\C[[\hbar]]$-module $E$ we
denote by $E_0$ its quotient $E/\hbar E$.

By a deformation of a complex vector space $E_0$ we mean a free $\C[[h]]$-module
$E$ such that $E/\hbar E \simeq E_0 $.
Deformation of an associative algebra
$\A_0$
is a $\C[[h]]$-algebra $\A$ such that $\A /\hbar\A \simeq \A_0 $
as associative $\C$-algebras.
The term {\em deformation quantization} or simply {\em quantization} is reserved for deformation of commutative algebras.

Quantized universal enveloping algebras are understood as
$\C[[\hbar]]$-algebras, \cite{Dr1}.
We will work with the standard or
Drinfeld-Jimbo quantization of simple Lie algebras.
We will assume $\hbar$-adic completion of the Cartan
subalgebra only. That is possible thanks to the existence
of the Poincar\'{e}-Birkhoff-Witt base over the Cartan subalgebra.
Thus defined, $\U_\hbar(\g)$ is a Hopf algebra in a weaker sense.
However that is sufficient
for our purposes, because we will deal with $\h$-diagonalizable
$\U_\hbar(\g)$-modules.

Although the quasitriangular structure requires completion
of tensor products, we will use this structure in
the situation when one of the modules is $\C[[\hbar]]$-finite.
This does not lead out of $\U_\hbar(\g)$.

We assume that $\U_\hbar(\g)$-modules are equipped with $\hbar$-adic topology
and the action of $\U_\hbar(\g)$ is continuous. However we do not require the modules
to be complete. We will work with $\U_\hbar(\g)$-modules that are direct sums
of $\C[[\hbar]]$-finite weight spaces. Each weight space is complete, being
$\C[[\hbar]]$-finite. Thus our point of view is self-consistent, because we
assume $\hbar$-adic completion for the Cartan subalgebra only.

\subsection{Quantized universal enveloping algebra}
\label{ssecQUEA}
Let $\g$ be a complex semisimple Lie algebra and $\h$ its Cartan subalgebra.
Let $\Rm$ denote the root system of $\g$ with a fixed subsystem of
positive roots $\Rm_+\subset \Rm$. By  $\Pi\subset \Rm_+$ we denote
the subset of simple roots.
Let $(.,.)$ denote the Killing form on $\g$. We will use the same notation for
the invariant scalar product on $\h^*$ that is induced by the Killing form restricted
to $\h$. By $\Wm$ we denote the Weyl group of $\g$.

For every $\la\in \h^*$ we denote by $h_\la$ its image under the isomorphism $\h^*\simeq \h$
implemented by the Killing form. In other words, $\la(h)=(h_\la,h)$ for all $h\in \h$.
We put $\rho=\frac{1}{2}\sum_{\al\in \Rm_+}\al $, the half-sum of the positive roots.
Its dual $h_\rho\in \h$
is the unique solution of the  system of linear equations
$\al_i(h_\rho)=\frac{1}{2}(\al_i,\al_i)$, $i=1,\ldots ,\rk\: \g$.

Denote by $\U_\hbar(\h)$ the $\hbar$-adic completion of the algebra $\U(\h)\tp \C[[\hbar]]$.
Define $\U_\hbar(\g)$ as a $\C[[h]]$-algebra generated by the
elements $e_{\pm\al}$, $\al\in \Pi $, over $\U_\hbar(\h)$. These generators are  subject
to the following relations, \cite{Dr1,Ji}:
$$
[h_{\al_i},e_{\pm\al_j}]=\pm(\al_i,\al_j) e_{\pm\al_j},
\quad
[e_{\al_i},e_{-\al_j}]=\delta_{ij}\frac{q^{h_{\al_i}}-q^{-h_{\al_i}}}{q_i-q^{-1}_i},
$$
$$
\sum_{k=0}^{1-a_{ij}}(-1)^k
\left[
\begin{array}{cc}
1-a_{ij} \\
 k
\end{array}
\right]_{q_i}
e_{\pm \al_i}^{1-a_{ij}-k}
e_{\pm \al_j}e_{\pm \al_i}^{k}
=0
,
$$
where  $a_{ij}=\frac{2(\al_i,\al_j)}{(\al_i,\al_i)}$,
$i,j=1,\ldots, \rk\: \g$, is the Cartan matrix, $q:=e^\hbar$, $q_i:= e^{\frac{\hbar}{2}(\al_i,\al_i)}$, and
$$
\left[
\begin{array}{cc}
n  \\ k
\end{array}
\right]_{q}
=
\frac{[n]_q!}{[k]_q![n-k]_q!},
\quad
[n]_q!=[1]_q\cdot [2]_q\ldots [n]_q
,
\quad
[n]_q=\frac{q^n-q^{-n}}{q-q^{-1}}.
$$

Strictly speaking, $\U_\hbar(\g)$ is not a Hopf algebra, because $\Delta$ takes values in
the completed tensor square of $\U_\hbar(\g)$. However, the category of $\U_\hbar(\g)$-modules
that are $\h$-diagonalizable is monoidal. So we will understand Hopf algebras in this weaker sense.

The coproduct $\Delta$, antipode $\gm$, and counit $\eps$ in $\U_\hbar(\g)$ are given by
\be
&\Delta(e_{\al})=e_{\al}\tp 1 + q^{h_{\al}}\tp e_{\al},
\quad
\Delta(e_{-\al})=e_{-\al}\tp q^{-h_{\al}} + 1 \tp e_{-\al},
\nn\\&
\Delta(h_{\al})=h_{\al}\tp 1 + 1\tp h_{\al},
\nn\\&
\gm(e_{\al})=-e_{\al}q^{-h_{\al}}
,\quad
\gm(e_{-\al})=-q^{h_{\al}}e_{-\al}
,\quad
\gm(h_{\al})=-h_{\al}.
\nn\\
&\eps(h)=\eps(e_{\pm \al})=0,
\nn
\ee
for all $\al \in \Pi$.
The correspondence $e_{-\al}\mapsto e_{\al}q^{h_\al}, e_{\al}\mapsto q^{-h_\al}e_{-\al}, h_\al\to h_\al$,
for all $\al\in \Pi$,
extends to an involutive anti-algebra and coalgebra transformation of $\U_\hbar(\g)$ denoted by $\omega$.

The $\C[[\hbar]]$-adic completion of $\U_\hbar(\g)$ is isomorphic to $\U(\g)[[\hbar]]$
as an associative algebra, \cite{Dr3}.

The subalgebras in $\U_\hbar(\g)$ generated over $\U_\hbar(\h)$ by $\{e_{+\al}\}_{\al\in \Pi} $ and
by $\{e_{-\al}\}_{\al\in \Pi} $,
respectively, are Hopf algebras.
They are quantized universal enveloping algebras of the positive and negative
Borel subalgebras $\b^\pm$ and denoted further by $\U_\hbar(\b^\pm)$.

The elements $\{e_{\pm\al}\}_{\al\in \Pi }\subset \U_\hbar(\g)$ are called quantum Chevalley  generators.
The Chevalley generators can be extended to a system of quantum Cartan-Weyl generators
$\{e_{\pm\al}\}_{\al\in \Rm_+}$ via the so called $q$-commutators, see \cite{KhT} and
references therein. The Cartan-Weyl base admits
an  ordering, with respect to which the monomials $\prod_{\al\in \Rm_+}^>e_{-\al}^{r_\al}\;\prod_{\al\in \Rm_+}^<e_{\al}^{s_\al}$
form a Poincar\'{e}-Birkhoff-Witt (PBW) base of $\U_\hbar(\g)$ as a $\U_\hbar(\h)$-module, \cite{KhT}.
Denote by $\U_\hbar'(\b^\pm)$ the ideals in $\U_\hbar(\b^\pm)$ generated by $\{e_{\pm\al}\}_{\al\in \Pi}$.
The PBW monomials of positive and, respectively, negative weights with respect to
$\h$ form bases for $\U_\hbar'(\b^\pm)$ over $\U_\hbar(\h)$.

The universal R-matrix belongs to the completed tensor square of $\U_\hbar(\g)$.
It has the structure
\be
\Ru=q^{\Omega_\h} \hspace{-7pt}\mod \U'_\hbar(\b^-)\hat\tp \U'_\hbar(\b^+),
\label{Rmat}
\ee
where $\Omega_\h\in \h\tp \h$ is the inverse to the Killing form (the canonical element)
restricted to $\h$. More precisely, the R-matrix can be represented as the product of the Cartan
factor $q^{\Omega_\h}$ and a series in the Cartan-Weyl generators, \cite{KhT}.
The universal R-matrix (\ref{Rmat}) is a "quantization" of the classical r-matrix
\be
\label{DJ-rm}
r_-=\frac{1}{2}\sum_{\al\in \Rm^+}  (\al,\al)(e_{-\al}\tp e_{\al}-e_{\al}\tp e_{-\al})
\ee
called the standard or Drinfeld-Jimbo r-matrix. Here $e_{\al}\in \g$ are the root vectors
normalized to $(e_{-\al},e_\al)=\frac{2}{(\al,\al)}$ with respect to the Killing form.
For $\al\in \Pi$, they are the classical limits of the quantum Chevalley  generators.

\subsection{Defining representations of classical matrix groups}
\label{ssecVR}
By classical Lie algebras (resp. algebraic groups) we mean the simple complex Lie algebras
of the types $\g=sl(n)$, $so(2n+1)$, $sp(n)$,  and $so(2n)$, for $n>1$.
We reserve the notation $V_0$ for the simple (defining) $\g$-module of dimension
$N=n$, $2n+1$, $2n$, and $2n$, respectively.
We chose the following realization of orthogonal algebras and symplectic algebras.
The algebra $so(N)$ leaves invariant
the skew-diagonal unit matrix; the symplectic algebra preserves the
skew-diagonal matrix with $+1$ above the center and $-1$ below.
We assume that $\h$, $\n^+$, and  $\n^-$ are realized in $\End(V_0)$ by, respectively,
diagonal, upper- and lower triangular matrices.

Bellow we collect some facts about the defining representation which we will use in
our exposition. Let $e_{ij}$ is the standard matrix base in $\End(V_0)$.
Define linear functionals $\{\ve_i\}_{i=1}^n$ on $\h$ setting
$\ve_i(e_{ii})=1$ for $\g=sl(n)$ and  $\ve_i(e_{ii}-e_{i'i'})=1$ otherwise.
Here $i'=N+1-i$.

The sets of positive and simple roots of $\g$
are expressed through $\ve_i$  by
$$
\begin{array}{lrlccl}
\Rm_+=&\{\ve_i-\ve_j\}_{i<j},&\Pi =\{\ve_1-\ve_2,\ldots, \ve_{n-1}-\ve_n\},& \g&=&sl(n),\\
\Rm_+=&\{\ve_i\pm\ve_j,\ve_i\}_{i<j},&\Pi =\{\ve_1-\ve_2,\ldots, \ve_{n-1}-\ve_n, \ve_n\},& \g&=&so(2n+1),\\
\Rm_+=&\{\ve_i\pm\ve_j,2\ve_i\}_{i<j},&\Pi =\{\ve_1-\ve_2,\ldots, \ve_{n-1}-\ve_n, 2\ve_n\},& \g&=&sp(n),\\
\Rm_+=&\{\ve_i\pm\ve_j\}_{i<j},&\Pi =\{\ve_1-\ve_2,\ldots, \ve_{n-1}-\ve_n,\ve_{n-1}+\ve_n\},& \g&=&so(2n).
\end{array}
$$
Our choice of functionals $\{\ve_i\}$ coincides with \cite{VO} excepting $\g=sl(n)$.
In this latter case $\{\ve_i\}$ are linearly dependent and satisfy
the condition $\sum_{i=1}^n \ve_i=0$.

The half-sum of positive roots is expressed in terms of $\{\ve_i\}$ by
$$
\rho=\sum_{i=1}^n \rho_i\ve_i,
\quad
\rho_i=\rho_1-(i-1),
\quad
\rho_1=
\left\{
\begin{array}{rcccl}
\frac{n-1}{2}&\mbox{for}& \g&=&sl(n),\\
n-\frac{1}{2}&\mbox{for}& \g&=&so(2n+1),\\
n&\mbox{for}& \g&=&sp(n),\\
n-1&\mbox{for}& \g&=&so(2n).\\
\end{array}
\right.
$$
The set of weights of the defining representation is
$\{\ve_i\}_{i=1}^n$ for $sl(n)$, $\{0\}\cup\{\pm\ve_i\}_{i=1}^n$ for $so(2n+1)$,
and $\{\pm\ve_i\}_{i=1}^n$ for $sp(n)$ and $so(2n)$.
\subsection{Parabolic subalgebras in $\U_\hbar(\g)$}
\label{ssecPS}
An element $\xi\in \g$ is called semisimple if $\ad\: \xi$ is diagonalizable.
A semisimple element belongs to a Cartan subalgebra, and all the Cartan subalgebras
are conjugated; so one can assume that $\xi\in \h$.

Let $\g=\n^-\oplus \h\oplus \n^+$ be the triangular decomposition
relative to $\Rm_+$.
A Levi subalgebra in $\g$ is defined as the centralizer of a semisimple element.
It is a reductive Lie algebra of the rank $\rk\; \g$.
Let $\l$ be a Levi subalgebra in $\g$ and $\p^\pm=\l +\n^\pm\subset \g$ be
the parabolic subalgebras.
Denote by $\n^\pm_\l$ the nillradicals in $\p^\pm$. The triangular decomposition
$\g=\n^-_\l \oplus \l \oplus \n^+_\l$ induces decomposition
$\U(\g)=\U(\n^-_\l)\U(\l)\U(\n^+_\l)$, which has a quantum analog, \cite{JT}.

The elements  $\{e_{\al},e_{-\al}\}_{\al\in \Pi_\l}$ generate over $\U_\hbar(\h)$
a Hopf subalgebra $\U_\hbar (\l)$ in $\U_\hbar (\g)$.
This subalgebra is a quantized universal enveloping algebra of the Levi subalgebra $\l\subset\g$.
It can be represented as  $\U_\hbar(\l_0)\U_\hbar(\c)$, where
$\l_0=[\l,\l]$ is the semisimple
part of $\l$ and  $\c\subset \l$ is the center.
Also, $\U_\hbar (\p^\pm):=\U_\hbar (\l)\U_\hbar (\b^\pm)$ are Hopf subalgebras.
They are quantized universal enveloping
algebras of the parabolic subalgebras in $\g$. This fact follows from the existence of
the PBW $\U_\hbar(\h)$-base.

Let $\Z_+$ denote the set of non-negative integers.
Consider in $\U_\hbar(\b^+)$  the sum of weight spaces with weights from
$\Z_+ (\Pi_\g-\Pi_\l)$. It is an algebra and
a deformation of $\U(\h+\n^+_\l)$, due to the existence
of the PBW base  in  $\U(\h+\n^+_\l)$. Let us denote this algebra by $\U_\hbar(\h+\n^+_\l)$.
According to \cite{Ke} (see also \cite{JT}), there is a subalgebra in $\U_\hbar(\h+\n^+_\l)$,
denoted further by $\U_\hbar(\n^+_\l)$, such that  $\U_\hbar(\b^+_\l)\U_\hbar(\n^+_\l)=\U_\hbar(\b^+)$.
Here $\b^+_\l=\l\cap \b_+$ is the positive Borel subalgebra in $\l$.
The algebra $\U_\hbar(\n^+_\l)$ is $\U_\hbar(\l)$-invariant
with respect to the adjoint action, and there exists a smash product decomposition
\be
\label{smash-}
\U_\hbar(\p^+)=\U_\hbar(\l)\ltimes \U_\hbar(\n^+_\l).
\ee
\begin{propn}
The algebra $\U_\hbar(\n^+_\l)$ is a deformation of $\U(\n^+_\l)$.
\end{propn}
\begin{proof}
Decomposition (\ref{smash-}) induces the decomposition
$\U_\hbar(\h+\n^+_\l)=\U_\hbar(\h)\ltimes \U_\hbar(\n^+_\l)\subset \U_\hbar(\p)$.
Since $\U_\hbar(\h+\n^+_\l)$ is a free $\U_\hbar(\h)$-module generated by the PBW base,
$\U_\hbar(\n^+_\l)$ is isomorphic to $\U(\n^+_\l)\tp \C[[\hbar]]$ as a $\C[[\hbar]]$-module.
By construction, the algebra $\U_\hbar(\n^+_\l)$ is generated by
$\ad(u)(e_{\al})$, where $\al\in \Pi_\g-\Pi_\l$ and $u\in  \U_\hbar(\l)$.
This implies the proposition.
\end{proof}

Denote by $\U_\hbar(\n^-_\l)$ the  $\omega$-image of $\U_\hbar(\n^+_\l)$, where
$\omega$ is the quantum Chevalley involution.
Since $\omega\bigl(\U_\hbar(\p^+)\bigr)=\U_\hbar(\p^-)$
and $\U_\hbar(\l)$ is $\omega$-stable, (\ref{smash-}) induces decomposition
$\U_\hbar(\p^-)=\U_\hbar(\n^-_\l)\rtimes \U_\hbar(\l)$, through $\omega$.
The algebra $\U_\hbar(\g)$ admits the triangular decomposition
\be
\U_\hbar(\g)=\U_\hbar(\n^-_\l)\U_\hbar(\l)\U_\hbar(\n^+_\l),
\label{q-triang}
\ee
which is a deformation of the classical one,  \cite{JT}.
Note that $\U_\hbar(\n^\pm_\l)$ are not Hopf algebras.

Thanks to the triangular decomposition (\ref{q-triang}), the algebra $\U_\hbar(\g)$ can be represented
as a direct sum
$$
\U_\hbar(\l)\oplus\sum_{\al\in \Pi_\g-\Pi_\l}(e_{-\al}\U_\hbar(\g)+\U_\hbar(\g)e_{\al}).
$$
By $\Pc_\l\colon \U_\hbar(\g)\to \U_\hbar(\l)$ we denote the projection along the second summand.

\section{Generalized Verma modules over $\U_\hbar(\g)$}
\label{secGVM}
In the present section we study tensor product
of finite dimensional and generalized Verma modules  over $\U_\hbar(\g)$.
This will be the basis for our further considerations.

\subsection{Upper and lower (generalized) Verma modules}
Let $\U_q(\g)$ be the quantum group in the sense of  Lusztig, \cite{L}.
It is a $\C(q)$-Hopf algebra generated by $\{e_{\pm\al_i},q_i^{\pm h_{\al_i}}\}_{\al_i\in \Pi}$.
The algebra $\U_q(\g)$ contains a $\C[q,q^{-1}]$-Hopf subalgebra $\breve\U_q(\g)$
such that $\U_q(\g)\simeq \breve\U_q(\g)\tp _{\C[q,q^{-1}]}\C(q)$. It is
generated by $\{e_{\pm\al_i}, q_i^{\pm h_{\al_i}}, [e_{\al_i},e_{-\al_i}]\}_{\al_i\in \Pi}$,
 see e.g. \cite{DCK}.
Clearly the algebra $\breve\U_\hbar (\g) :=\breve \U_q(\g)\tp_{\C[q,q^{-1}]}\C[[\hbar]]$,
where $\C[q,q^{-1}]$ is embedded in $\C[[\hbar]]$ via $q\mapsto e^\hbar$, is
dense in $\U_\hbar(\g)$ in the $\hbar$-adic topology. Remark that $\breve \U_\hbar(\g)$
is $\U_\hbar(\g)$-invariant with respect to the adjoint action.

There is a one-to-one correspondence between  finite-dimensional $\g$-modules and finite-dimensional $\U_q(\g)$-modules
with $q^{\Z}$-valued weights. Each such module is isomorphic to $\breve W\tp_{\C[q,q^{-1}]}\C(q)$,
where $\breve W$ is a $\breve \U_q(\g)$-module, free and finite over $\C[q,q^{-1}]$, \cite{Jan2}.
The specialization $\breve W \mod (q-1)$ gives a finite dimensional $\g$-module.
Therefore $\breve W$ extends to a $\U_\hbar(\g)$-module, free and finite over $\C[[\hbar]]$.
We will call such $\U_\hbar(\g)$-modules \select{finite dimensional}. They are deformations of $\g$-modules,
diagonalizable over $\U_\hbar(\h)$, and have the same weight
structure. This correspondence between finite dimensional  $\g$-modules and $\U_\hbar(\g)$-modules is additive.
That is, every finite dimensional $\U_\hbar(\g)$-module is a direct sum of "simple" modules,
i. e. those whose classical limit is simple.

Every submodule of a finite dimensional $\U_\hbar(\g)$-module is $\h$-diagonalizable
with $q^{\Z}$-valued weights; hence it is again finite dimensional.
We will also deal with $\U_\hbar(\g)$-modules that are finitely generated over $\C[[\hbar]]$ but
not free. We always assume that such modules are quotients of finite dimensional
and called them just $\C[[\hbar]]$-finite.

A highest  weight $\U_\hbar(\g)$-module is generated by a weight vector
annihilated by $\U_\hbar'(\n^+)$. Similarly, a lowest weight module
is generated by a weight vector annihilated by $\U_\hbar'(\n^-)$.
Finite dimensional $\U_\hbar(\g)$-modules
have highest and lowest weights simultaneously.
Their highest weights are integral dominant.
Finite dimensional $\U_\hbar(\g)$-modules
with highest weights are deformations of irreducible finite dimensional $\g$-module.
They are substitutes for irreducibles, since a $\U_\hbar(\g)$-module
is almost never irreducible in the usual sense (multiplication by $\hbar$ is a morphism).

For reductive $\l$, the highest weight $\U_\hbar (\l)$-modules are defined similarly to
the semisimple case.
Let $A$ be a representation of  $\U_\hbar (\l)$
with highest weight.
It extends to a representation of $\U_\hbar (\p)$ for $\p=\p^+$,
by setting it zero on $\U_\hbar' (\n^+_\l)$.
A generalized Verma module over $\U_\hbar (\g)$ is the induced module
$M_{\p,A}=\U_\hbar (\g)\tp_{\U_\hbar (\p)} A=:\Ind^\g_\p\: A$, cf. \cite{J1}.
By $\C_\la$ we will denote the one dimensional $\U_\hbar(\l)$-module defined by a
character $\la\colon \l\to \C$.
The module $M_{\p,\C_\la}$ will be denoted simply by $M_{\p,\la}$ and
the Verma module  $M_{\b,\la}$ by $M_{\la}$.

Every $\U_{\hbar}(\g)$-module with highest weight $\la$ is a quotient of
$M_{\la}$. The generalized Verma module $M_{\p,A}$ induced from a
$\U_\hbar(\l)$-module  $A$ with highest weight $\la$ is a highest weight module, so it is
a quotient of $M_{\la}$ as well.

\begin{lemma}
\label{cl-q}
Let $A$ be a finite dimensional $\U_\hbar(\l)$-module. Then the
$\U_\hbar(\g)$-module $M_{\p,A}$ is a deformation of the
classical generalized Verma module over $\U(\g)$.
\end{lemma}
\begin{proof}
It follows from (\ref{q-triang}) that $M_{\p,A}$
is a free $\U_\hbar(\n^-_\l)$-module generated
by $1\tp_{\U_\hbar (\p)} A$ and hence it is $\C[[\hbar]]$-free.
Since the decomposition (\ref{q-triang}) is a deformation of the classical triangular decomposition,
$M_{\p,A}$ is a deformation of the corresponding generalized Verma module over
$\U(\g)$.
\end{proof}

A highest weight representation of $\U_\hbar (\l)$ can be extended to a representation
of $\U_\hbar (\p^-)$ by setting it trivial on $\U_\hbar(\n_\l^-)$.
This is possible, due to (\ref{smash-}).
Similarly to
$M_{\p^+,A}$, the  module $M_{\p^-,A}:=\U_\hbar(\g)\tp _{\U_\hbar(\p^-)} A$
is introduced, where $A$ is taken to be a lowest weight $\U_\hbar(\l)$-module.
\begin{propn}
\label{+-l}
Let $A$ be a lowest and $B$ a highest weight $\U_\hbar(\l)$-modules.
Then the $\U_\hbar(\g)$-module $M_{p^-,A}\tp M_{p^+,B}$ is
isomorphic to $\Ind^{\g}_{\l}(A\tp B):=\U_\hbar(\g)\tp_{\U_\hbar(\l)} (A\tp B) $.
\end{propn}
\begin{proof}
Consider the map $\U_\hbar(\g)\tp A\tp B\mapsto \bigl(\U_\hbar(\g)\tp A\bigr)\tp \bigl(\U_\hbar(\g)\tp B\bigr)$
defined by
$u\tp a\tp b\mapsto u^{(1)}\tp a\tp u^{(2)}\tp b$,
where $u^{(1)}\tp u^{(2)}$ is the standard symbolic notation for the coproduct $\Delta(u)$.
This map induces a homomorphism
\be
\label{F}
\Delta^\g_{\l,A,B}\colon\U_\hbar(\g)\tp_{\U_\hbar(\l)} (A\tp B) \to M_{p^-,A}\tp M_{p^+,B}
\ee
 of $\U_\hbar(\g)$-modules. We claim that this map is an isomorphism.

First let us prove the statement assuming $\l=\h$, $A=\C_\mu$, and $B=\C_\nu$.
Introduce the grading in $\U_\hbar(\n^\pm)$ by weight height setting $\deg e_{\pm\al}=1$ for $\al\in \Pi$.
The gradings in $\U_\hbar(\n^\pm)$ induce a double grading in $\Ind^\g_\h\bigl(\C_\mu \tp \C_\nu\bigr)$
and $M_{p^-,\mu}\tp M_{p^+,\nu}$, which
can be identified as graded spaces. Each homogeneous component has finite rank over $\C[[\hbar]]$.
Let $v_\mu$ and $v_\nu$ be the generators of  $M_{p^-,\mu}$
and $M_{p^+,\nu}$, respectively.
Take $u^\pm\in \U_\hbar(\n^\pm)$ to be monomials in $\{e_{\pm\al}\}_{\al\in \Pi}$
and compute the map $\Delta^\g_{\l,\C_\mu,\C_\nu}$:
\be
 u^-u^+(v_\mu\tp v_\nu)&\mapsto& (u^-)^{(1)}(u^+)^{(1)}v_\mu\tp (u^-)^{(2)}(u^+)^{(2)}v_\nu=
(u^-)^{(1)}u^+ v_\mu\tp (u^-)^{(2)} v_\nu
\nn\\
&=&c_q u^+v_\mu\tp u^-v_\nu + w^+v_\mu\tp w^-v_\nu.
\nn
\ee
Here $c_q\in \C[[\hbar]]$ is invertible and the elements $w^\pm \in \U_\hbar(\n^\pm)$ belong
to subspaces of degree $<\deg u^\pm$.
This computation shows that $\Delta^\g_{\l,\C_\mu,\C_\nu}$ is a triangular operator (relative to the double grading)
with invertible diagonal. Therefore it is an isomorphism.

The above consideration also proves that the map (\ref{F})  is an epimorphism in the general situation,
as $A$ and $B$ are quotients of
$M_{p^-,\mu}$ and $M_{p^+,\nu}$, respectively.
We must check that (\ref{F}) is injective.
The map $\Delta^\g_{\l,A,B}$ is surjective modulo $\hbar$ as a $\U(\g)$-morphism.
By dimensional arguments based on the bi-grading, we conclude that this $\U(\g)$-morphism is an isomorphism.
Therefore $\Delta^\g_{\l,A,B}$  is an isomorphism by the obvious deformation arguments.
\end{proof}
\begin{remark}
\label{breve}
In conclusion of this section we remark that the triangular decomposition and the generalized
Verma modules can by naturally defined for the algebra $\breve \U_\hbar(\g)$. We will use
this observation in Section \ref{secQCC}, where the induction is made from the character $\la/2\hbar$,
$\la\in \c^*$.
The $\breve \U_\hbar(\g)$-action on $M_{\p,\la/2\hbar}$
does not extend to an action of $\U_\hbar(\g)$.
 However the $\breve \U_\hbar(\g)$-action on $\End(M_{\p,\la/2\hbar})$ extends
to an action of $\U_\hbar(\h)$, and that is what we need for our construction.
\end{remark}
\subsection{Pairing between upper and lower generalized Verma modules}
Denote by $\c$ the center of $\l$ and by $\c_{reg}\subset \c$ the subset of elements whose centralizer
is exactly $\l$. Clearly $\c_{reg}$ is a dense open set in $\c$. The coadjoint $\g$-module $\g^*$ is canonically identified
with $\g$ via the Killing form.
Then the dual space $\c^*$ is identified with
the orthogonal complement to the annihilator of $\c$ in $\h^*$, so $\c^*\subset \h^*$ under this convention.

Given a root $\al\in \Rm$, let $\al^\vee$ denote the  dual root $\al^\vee=\frac{2}{(\al,\al)}\al$.
We call a weight $\la\in \c^*$ generic if
$(\la,\al^\vee)\not \in \Z$ for all $\al\in \Pi_\g-\Pi_\l$. Clearly the set $\c^*_{gen}$ of
generic weights is a dense open subset in
$\c^*_{reg}=\{\la\in \c^*\:|\: (\la,\al^\vee)\not =0, \forall \al\in \Pi_\g-\Pi_\l\}$.

With every finite dimensional irreducible  $\l$-module $A$ one can associate
a weight $\la_A\in \c^*$ such that $h a=\la_A(h)a$ for all $h\in \c$ and $a\in A$.
We call $A$ generic if $\la_A\in\c^*_{gen}$. An arbitrary finite dimensional
$\l$-module is called generic if its every irreducible submodule is generic.
This terminology extends to the corresponding $\U_\hbar(\l)$-modules.

There exists a $\U_\hbar(\g)$-equivariant pairing between $M_{p^-,A^*}$ and $M_{p^+,A}$.
The construction goes as follows.
Consider a bilinear $\U_\hbar(\g)$-equivariant map $(\U_\hbar(\g)\tp A^*)\tp (\U_\hbar(\g)\tp A)\to \C[[\hbar]]$ defined
by $u_1\tp \xi \tp u_2\tp x\mapsto \xi\Bigl(\Pc_\l\bigl(\gm(u_1)u_2\bigr)x\Bigr)$,
where $\Pc_\l$ is introduced in Subsection \ref{ssecPS}.
This map is equivariant by construction and factors through a bilinear equivariant map
$M_{p^-,A^*}\tp M_{p^+,A}\to \C[[\hbar]]$,
as required.

\begin{propn}
\label{Shap}
Let $A$ be a generic finite dimensional $\U_\hbar(\l)$-module.
Then the equivariant pairing between
$M_{\p^+,A}$ and $M_{\p^-,A^*}$
is nondegenerate.
\end{propn}
\begin{proof}
Without loss of generality, we may assume that $A_0=A/\hbar A$ is irreducible.

Since $-\c^*_{gen}=\c^*_{gen}$, the modules $A$ and $A^*$ are generic simultaneously.
It follows from \cite{Jan1}, Satz 3, that the $\U(\l)$-module $M_{\p^+,A}/\hbar M_{\p^+,A}$
is irreducible  for $\la\in \c^*_{gen}$. Clearly the same is
true for $M_{\p^-,A^*}/\hbar M_{\p^-,A^*}$.
Therefore the pairing in question
is nondegenerate modulo $\hbar$ (being $\U(\l)$-equivariant and not  identically zero).

With respect to the pairing, the weight spaces of weights $\mu$ and $\nu$ are orthogonal
unless $\mu+\nu=0$. Since the weight spaces in $M_{\p^+,A}$ and $M_{\p^-,A^*}$
are $\C[[\hbar]]$-finite, non-degeneracy of the pairing follows from
non-degeneracy modulo $\hbar$.
\end{proof}

\subsection{Tensor product of finite dimensional and generalized Verma modules}
We call a $\U_\hbar(\g)$-module a weight module if it is
$\h$-diagonalizable and its weight spaces are finite and free over $\C[[\hbar]]$.

For any  weight module let $\Lambda'(U)$ denote the set
of weights of $U$ and $U[\mu]$ the weight space for $\mu\in \Lambda'(U)$.
The dual module $U^*$ of linear functionals on $U$
consists of infinite formal sums $f=\sum_\mu f_\mu$, where $f_\mu\in U[\mu]^*$.
The action of $\U_\hbar(\g)$ on $U^*$ is defined to be $(x f)(u)=f(\gm(x)u)$,
for $f\in U^*$, $u\in U$, and $x\in \U_\hbar(\g)$.
Define the restricted dual $U^\circ$ as a natural $\U_\hbar(\g)$-submodule in $U^*$
by setting  $U^\circ = \oplus_{\mu\in \Lambda'(U)} U[\mu]^*\subset U^*$ (only finite sums admitted).

\begin{lemma}
\label{rdual}
Let $U_1,U_2$ be  weight $\U_\hbar(\g)$-modules  and $W_1,W_2$ finite dimensional $\U_\hbar(\g)$-modules. Then
$\Hom_{\U_\hbar(\g)}(W_1\tp U_1, W_2\tp U_2)\simeq \Hom_{\U_\hbar(\g)}(U_2^\circ\tp U_1, W_2\tp W_1^*)$.
\end{lemma}
\begin{proof}
Clear.
\end{proof}
Proposition \ref{Shap} asserts that $M_{\p^-,A^*}\simeq M^\circ_{\p,A}$ for generic $A$.
\begin{lemma}
\label{aux}
Let $W$ be a finite dimensional $\U_\hbar(\g)$-module,  $A$ and $B$ finite dimensional $\U_\hbar(\l)$-modules,
and $\la\in \c^*_{gen}$.
Then the following $\C[[\hbar]]$-linear isomorphisms  take place:
\be
\label{Hom_isom}
\begin{array}{lcl}
\Hom_{\U_\hbar(\g)}(M_{\p,A\tp \C_\la},M_{\p,B\tp \C_\la})&\simeq &\Hom_{\U_\hbar(\l)}(A,B),
\\
\Hom_{\U_\hbar(\g)}(M_{\p,A\tp \C_\la},W\tp M_{\p,\la})&\simeq &\Hom_{\U_\hbar(\l)}(A,W),
\\
\Hom_{\U_\hbar(\g)}(W\tp M_{\p,\la},M_{\p,A\tp \C_\la})&\simeq &\Hom_{\U_\hbar(\l)}(W,A),
\\
\Hom_{\U_\hbar(\g)}(W\tp M_{\p,\la},W\tp M_{\p,\la})&\simeq &\Hom_{\U_\hbar(\l)}(W,W).
\end{array}
\ee
\end{lemma}
\begin{proof}
The proof is based on Propositions \ref{+-l} and  \ref{Shap} and can be conducted similarly as
in \cite{DM1} for the classical case of $\U(\g)$-modules.
For instance, let us check the first isomorphism. By Proposition \ref{Shap}, the module
$M_{\p^-,\C^*_\la\tp B^*}$ is isomorphic to the restricted dual
$M_{\p,B\tp \C_\la}^\circ$.
Therefore
\be
\label{eq_aux}
\Hom_{\U_\hbar(\g)}(M_{\p,A\tp \C_\la},M_{\p,B\tp \C_\la})\simeq
\Hom_{\U_\hbar(\g)}(M_{\p^-,\C^*_\la\tp B^*}\tp M_{\p^+,A\tp \C_\la},\C[[\hbar]]),
\ee
 by Lemma \ref{rdual}.
Here $\C[[\hbar]]$ is the trivial $\U_\hbar(\g)$-module.
According to Lemma \ref{+-l},
the tensor product of lower and upper generalized Verma modules is
induced from the $\U_\hbar(\l)$-module $\C^*_\la\tp B^* \tp A\tp \C_\la\simeq B^*\tp A$.
Applying the Frobenius reciprocity, we continue
(\ref{eq_aux}) with $\Hom_{\U_\hbar(\l)}(B^*\tp A,\C[[\hbar]])\simeq \Hom_{\U_\hbar(\l)}(A,B)$, as required.
\end{proof}
\begin{propn}
\label{dir_sum}
Let $W$ be a finite dimensional $\U_\hbar(\g)$-module and $\la\in \c^*_{gen}$. Then
$W\tp M_{\p,\la}$ admits the direct sum decomposition
\be
\label{eq_dir_sum}
W\tp M_{\p,\la}= \oplus_{A_0} M_{\p,A\tp\C_\la},
\ee
where summation is taken over the simple $\l$-modules
with multiplicities entering $W_0$.
\end{propn}
\begin{proof}
The isomorphisms (\ref{Hom_isom}) hold modulo $\hbar$. Moreover,
they  commute with taking quotients $\mod \hbar$.
First let us prove  the classical $\mod \hbar$ analog of decomposition (\ref{eq_dir_sum}), retaining
the same notation for the Verma modules over $\U(\g)$.
We assume tn  (\ref{eq_dir_sum}) a fixed decomposition within
each isotypic $A_0$-component of $W_0$.
Let $j^{A_0}\colon A_0 \hookrightarrow W_0$ be the $\l$-equivariant injections
such that $\sum_A j^{A_0}=\id_{W_0}$. Let ${\hat j}^{A_0}$
be their lifts $M_{\p,A_0\tp \C_\la}\to W\tp M_{\p,A_0\tp \C_\la}$.
For generic $\la\in \c^*_{gen}$ all the modules $M_{\p,A_0\tp \C_\la}$ are irreducible.
Therefore all ${\hat j}^{A_0}$ are linearly independent, and the $\g$-equivariant map
$ \oplus_{A_0} M_{\p,A_0\tp\C_\la}\stackrel{\sum_{A_0} {\hat j}^{A_0} }{\longrightarrow} W_0\tp M_{\p,\la}$ is an embedding.
Then  it is an embedding of $\h$-modules $\U(\n^-_\l)\tp W_0 \simeq \oplus_{A_0} \U(\n^-_\l)\tp A_0\stackrel{\sum_{A_0}
{\hat j}^{A_0} }{\longrightarrow} W_0\tp \U(\n^-_\l) $
and therefore an isomorphism. This proves the statement modulo $\hbar$.

We can chose $j^A\in \Hom_{\U_{\hbar}(\l)}(A,W)$ to be deformations of morphisms  $j^{A_0}$
splitting the $\U_{\hbar}(\l)$-module $W$ into the direct sum of highest weight submodules.
Take the $\U_\hbar(\g)$-morphisms ${\hat j}^A\colon M_{\p,A\tp \C_\la}\to W\tp M_{\p,\la}$
corresponding to $j^A\in \Hom_{\U_{\hbar}(\l)}(A,W)$ under the second isomorphism from (\ref{Hom_isom}).
Since isomorphisms (\ref{Hom_isom}) commute with taking quotients $\mod \hbar$,
we have ${\hat j}^A ={\hat j}^{A_0}\mod \hbar$.
Consider the morphism
$ \oplus_{A_0} M_{\p,A\tp\C_\la}\stackrel{\sum_{A_0} {\hat j}^A}{\longrightarrow} W\tp M_{\p,\la}$
of $\U_\hbar(\g)$-modules. Restricting consideration
to weight spaces we conclude that $\sum_{A_0} {\hat j}^A$ is an isomorphism because it is so modulo $\hbar$.
\end{proof}

\section{Properties of the universal RE matrix}
\label{secREM}
In the present section we recall general properties of the universal reflection equation (RE) matrix $\Q=\Ru_{21}\Ru$
and study its action on tensor products of finite
dimensional and (generalized) Verma modules.

Recall from \cite{Dr2} that the element $\Q=\Ru_{21}\Ru$ can be represented as
 $\Q=\Delta(\upsilon)(\upsilon^{-1}\tp \upsilon^{-1})$, where
$\upsilon=\gm(\Ru^{-1}_1) \Ru_2^{-1}=\gm^2(\Ru_1) \Ru_2$. Conjugation with
$\upsilon$ implements the squared antipode $\upsilon x\upsilon^{-1}=\gm^2(x)$ for all $x\in \U_\hbar(\g)$.
On the other hand, the squared antipode can be written as the conjugation
$\gm^2(x)= q^{-2h_\rho}x q^{2 h_\rho}$, $x\in \U_\hbar(\g)$.
It follows from here that
$\upsilon=q^{-2h_\rho}\zz$, where $\zz$ is some invertible element from the center
of $\U_\hbar(\g)$. Therefore we can write
\be
\Q=\Delta(\zz)(\zz^{-1}\tp \zz^{-1}),
\label{Q-Z}
\ee
as $q^{-2h_\rho}$ is a group-like element.
Let $\la\in \h^*$ and $\chi^\la$ be the corresponding central character of $\U_\hbar(\g)$.
It is easy to compute the value $\chi^\la(\zz)$ via a $\U_\hbar(\g)$-module
$W$ with highest weight, using the structure of R-matrix (\ref{Rmat}).
Let $w_\la$ be the highest weight vector in $W$.
Since $\U'_\hbar(\n^+)w_\la=0$, we have $\upsilon w_\la = q^{(\la,\la)}w_\la=\chi^\la(\zz) q^{-2(\rho,\la)}w_\la$.
From this we find
\be
\chi^\la(\zz)= q^{(\la,\la)+2(\rho,\la)}.
\label{la-Z}
\ee
The element $\Q$ satisfies the identities
\be
\Ru_{21}\Q_{13} \Ru_{12} \Q_{23}&=&\Q_{23}\Ru_{21} \Q_{13}\Ru_{12},
\label{ure}
\\
(\Delta\tp \id)(\Q)&=&\Ru^{-1}_{12}\Q_{13}\Ru_{12}\Q_{23},
\label{dQ1}
\ee
of which the first may be called the universal reflection equation,
cf. \cite{DKM}.
Equation (\ref{dQ1}) is the key identity of the fusion procedure for
solutions to the RE, \cite{DKM}.

\subsection{Characteristic polynomial for RE matrix}
Given a finite dimensional representation $(W,\pi_W)$ let $\Q_W$  denote the element
$(\pi_W\tp \id)(\Q)\in \End(W)\tp \U_\hbar(\g)$.
\begin{propn}
Let $W$ be a highest weight $\U_\hbar(\g)$-module with the
multiset of weights $\Lambda(W)$.
Then i) there exists a polynomial $p$ of degree $\# \Lambda(W)$ in one
variable with coefficients in the center of $\U_\hbar(\g)$
such that $p(\Q_W)=0$,
ii) the spectrum of the operator $\Q_{W}$
on $W\tp M_{\la}$
is
\be
\label{spec_V}
\bigl\{\:
q^{2(\la+\rho,\nu_i)-2(\rho,\nu)+(\nu_i,\nu_i)-(\nu,\nu)}
\:\bigr\}_{\nu_i\in \Lambda(W)},
\ee
where $\nu$ is the highest weight of $W$.
\label{spec}
\end{propn}
\begin{proof}
Notice that a symmetric function in the eigenvalues  (\ref{spec_V}) is invariant under the action
of the Weyl group.
Then i) follows from ii) through the Harish-Chandra homomorphism, \cite{Jan2}.
So let us check ii).

As a $\C[[\hbar]]$-module, $W\tp M_{\la}$ is isomorphic
to $W\tp \U_\hbar(\n^-)$. Under this isomorphism,
the subspace of weight $\bt\in \h^*$ in $W\tp \U_\hbar(\n^-)$
goes to the subspace of weight $\bt+\la$ in $W\tp M_{\la}$.
Since $\Q$ is invariant, it preserves the weight spaces
in $W\tp M_{\la}$, which have finite rank over $\C[[\hbar]]$.
Now observe that for generic $\la$ ii) follows from Proposition \ref{dir_sum}.
This implies ii) for arbitrary $\la$, since
restriction of $\Q_W$ to any weight space has polynomial dependance on $\la$ in every order in $\hbar$.
\end{proof}

\subsection{Minimal polynomial for RE matrix}
In this subsection we determine the spectrum of $\Q$ on the tensor product of finite dimensional and generalized
Verma modules over $\U_\hbar(\g)$.
\begin{thm}
\label{Spec_parabolic}
Let $\la\in \c^*_{gen}$ be  a generic character of a Levi subalgebra in $\l\subset \g$.
Let $W$ be a finite dimensional $\U_\hbar(\g)$-module and let
$\Lambda_\l(W)=\{\upsilon_l=\nu_{i_l}\}\subset \Lambda (W)$ denote the multiset of highest
weights of simple $\l$-modules entering $W_0$.
The operator $\Q$ is diagonalizable on $W\tp M_{\p,\la}$ and has
eigenvalues
\be
\label{roots_min}
\bigl\{\:q^{2(\la+\rho,\upsilon_l)-2(\rho,\nu)+(\upsilon_l,\upsilon_l)-(\nu,\nu)}\:\bigr\}_{\upsilon_l\in \Lambda_\l(W)}.
\ee
\end{thm}
\begin{proof}
For generic $\la$ the module $W\tp M_{\p,\la}$
splits into the direct sum of highest weight modules, by Proposition \ref{dir_sum}.
The operator $\Q$ is proportional to  $\Delta(z)$ on $W\tp M_{\p,\la}$.
Thus we conclude that  $\Q$ is diagonalizable because $\Delta(z)$  is diagonalizable.

For every $\U_\hbar(\l)$-module $A\subset W$ with highest weight $\mu$ the summand
$M_{\p,A\tp \C_\la}$ in (\ref{eq_dir_sum}) is a $\U_\hbar(\g)$-module of highest weight $\la+\mu$. Hence
$\Delta(z)$ and  therefore $\Q$ act as scalar multipliers on $M_{\p,A\tp \C_\la}$.
Their eigenvalues are computed using (\ref{la-Z}).
\end{proof}

\subsection{A construction of central elements}
Assume that $\g$ is an arbitrary complex simple Lie algebra.
Let $(W,\pi_W)$ be a finite dimensional $\U_\hbar(\g)$-module. Let $X\in \End(W)\tp \U_\hbar (\g)$
be an invariant matrix, i. e. commuting with
$\pi_W(x^{(1)})\tp x^{(2)}$ for all $x\in \U_\hbar (\g)$.
It is known that the q-trace defined by
$$\Tr_{q}(X):=\Tr\bigl(\pi_W(q^{2h_\rho})X\bigr)$$
is $\ad$-invariant and hence belongs to the center of $\U_\hbar (\g)$.

The annihilator of the Verma module with the highest weight $\la$ is
generated by the kernel of a central character $\chi^\la$. Let us compute the values
 $\chi^\la\bigl(\tau^\ell\bigr)$.

Define a map $d\colon \h^*\to \C[[\hbar]]$ setting
\be
d(\la):=\prod_{\al\in \Rm_+}\frac{ q^{(\la+\rho,\al)}-q^{-(\la+\rho,\al)}}{q^{(\la,\al)}-q^{-(\la,\al)}},
\quad q=e^\hbar.
\ee
For a finite dimensional $\U_\hbar(\g)$-module $W$ with the highest weight $\la$,
the Weyl character formula \cite{N} gives $d(\la)=\Tr_q(\id_W)$, the q-dimension of $W$.
\begin{propn}
Let $W$ be a
finite dimensional $\U_\hbar(\g)$-module with the multiset of weights $\Lambda(W)$.
Then for any $\la\in \h^*$
\be
\chi^\la\Bigl(\Tr_q(\Q_W^\ell)\Bigr)=
\sum_{\nu_i\in \Lambda(W)} x_{\nu_i}^\ell \frac{d(\la+\nu_i)}{d(\la)},
\label{char_V}
\ee
where $x_{\nu_i}$ are the eigenvalues of $\Q_W$ given by (\ref{spec_V}).
\end{propn}
\begin{proof}
We adapt a proof from \cite{GZB}, because it is suitable for
any $\U_\hbar(\g)$-invariant operator.

Let us check the statement for special $\la$ first. Namely,
we suppose that $\la$ is integral dominant and $\nu$, the highest weight of $W$, is subordinate to $\la$.
The weight $\nu$ is called subordinate to $\la$ if $\la+\nu_i$ are dominant integral
for all $\nu_i \in \Lambda(W)$. For $\nu$ fixed we denote  by $\D^+_\nu$ the set of such $\la$ that
$\nu$ is subordinate to $\la$.
According to  \cite{K}, a polynomial function on $\h^*$ is determined by its values on $\D^+_\nu$.
Both sides of (\ref{char_V}) are polynomials in $\la$ in  every order in $\hbar$,
thus it suffices to compare them on $\D^+_\nu$ only.

For $\la\in \D^+_\nu$ let us compute the central character in question
on the finite dimensional module $U$ with the highest weight $\la$.
In this case, the module $W\tp U$ splits into the direct sum
of $\U_\hbar(\g)$-modules of highest weights $\nu_i+\la$ for $\nu_i\in \Lambda(W)$.
An invariant operator  on $W\tp U$ decomposes into over invariant
projectors $\{P_i\}_{\nu_i\in \Lambda(W)}$ to the highest weight submodules, so it suffices to compute
the q-trace of these projectors.
The operator $(\Tr_q\tp \id)(P_i)$ is constant on $U$.
Taking q-trace over  $U$ gives
$(\Tr_q\tp \id)(P_i)=(\Tr_q\tp \Tr_q)(P_i)/(\id\tp \Tr_q)(P_i)=d(\la+\nu_i)/d(\la)$,
since the q-trace is multiplicative with respect to the tensor product.
This proves the statement for $\la\in \D^+_\nu$ and therefore for all $\la\in \h^*$.
\end{proof}

\section{On quantization of affine algebraic varieties}
\label{secQAAV}
In this section we develop a machinery for equivariant quantization.
Throughout the section we assume that $\g$ is semisimple.
Moreover, here we admit  an \select{arbitrary}, i. e.
triangular or even trivial, quantization of $\U(\g)$.

In what follows, we use some standard facts from commutative algebra, such
as the Nakayama lemma. The reader can find the "list" of facts we rely on in \cite{M}.

\subsection{A flatness criterion}
\label{secQn}
Recall that an associative algebra and (left) $\U_\hbar(\g)$-module $\A$ is called a $\U_\hbar(\g)$-module algebra if
its multiplication is compatible with the $\U_\hbar(\g)$-action.
That is, for all $h\in \U_\hbar(\g)$ and all $a,b\in \A$
$$
h\tr (ab)=(h^{(1)}\tr a)(h^{(2)}\tr b),
\quad\mbox{where}\quad h^{(1)}\tp h^{(2)}=\Delta(h).
$$
\begin{definition}
 A deformation $\A$ of a $ \U(\g)$-module algebra $\A_0$ is called equivariant if
$\A$ is a $ \U _\hbar(\g)$-module algebra and the action of
 $ \U _\hbar(\g)$ on $\A$ coincides modulo $\hbar$ with
 the  action of $ \U (\g)$ on $\A_0$.
\end{definition}

For every $\U_\hbar(\g)$-module $E$ and a finite
dimensional highest weight  $\U_\hbar(\g)$-module $W$ there exists a natural morphism
$W\tp \Hom_{\U_\hbar(\g)}(W,E)\to E$ of $\U_\hbar(\g)$-modules.
We call the image of this morphism the  isotypic
$W$-component of $E$.

We call a $\U_\hbar(\g)$-module  \select{admissible} if it is
a direct sum of its isotypic components and each of them is finitely generated
over $\C[[\hbar]]$. It can be shown that submodules and quotient
modules of an admissible module are admissible.

The following proposition gives a flatness criterion
for quotient algebras.
Informally, if one constructs
an ideal $\Jg$ that lies in the kernel of a certain homomorphism $\Sg\to \Tg$
and has the "right classical limit", then $\Jg$
equals the entire kernel.
\begin{propn}[deformation method]
\label{Lemma_aux}
Let $\Sg$ be an admissible and $\Tg$ torsion free $\U_\hbar(\g)$-module algebras;
let $\Psi\colon \Sg\to \Tg$ be a non-zero equivariant homomorphism.
Suppose $\ker \Psi$ contains an invariant ideal  $\Jg$ such that
the image $\Jg_0^\flat$ of $\Jg_0$  in $\Sg_0$ is a
maximal $\g$-invariant ideal.
Then  i) $\im \Psi$ is a $\C[[\hbar]]$-free deformation of $\Sg_0/\Jg_0^\flat$,
ii) $\ker \Psi=\Jg$.
\end{propn}
\begin{proof}
The $\U_\hbar(\g)$-module $\im \Psi$ is admissible and torsion free;
hence it is  $\C[[\hbar]]$-free.

Since $\im \Psi$ is free, we have the direct sum decomposition
$\Sg\simeq \ker \Psi\oplus \im \Psi$ of $\C[[\hbar]]$-modules.
Hence $\Sg_0\simeq (\ker \Psi)_0\oplus (\im \Psi)_0$
and $(\ker \Psi)_0\supset \Jg_0^\flat$. By assumption, $\Jg_0^\flat$ is a maximal $\g$-invariant
ideal, hence either $(\ker \Psi)_0 = \Sg_0$ or $(\ker \Psi)_0= \Jg_0^\flat$.
The first option is impossible. Indeed, then $(\im \Psi)_0=0$ and $\im \Psi=0$ since
$\im \Psi$ is $\C[[\hbar]]$-free.
In this case $\ker \Psi=\Sg$,  and the map $\Psi$ would have been zero.
Therefore $(\im \Psi)_0=\Sg_0/\Jg_0^\flat$ and i) is proven.

We have an embedding  $\Jg\hookrightarrow \ker \Psi$ and an epimorphism $\Jg_0\to (\ker \Psi)_0=\Jg_0^\flat$.
Applying the Nakayama lemma to each isotypic component, we prove ii).
\end{proof}

\begin{remark}
1. We emphasize that $\Jg_0^\flat$ is assumed to be not just $\Jg_0=\Jg/\hbar\Jg$ but its image in $\Sg_0$.
It is essential to distinguish between $\Jg_0$ and $\Jg_0^\flat$  because the functor $\mod \hbar$ is not left exact.
Eventually, $\Jg_0$ and $\Jg_0^\flat$ coincide in our situation. However that is not {\em a priory} obvious
and follows from the proof.
In practice, $\Jg$ is often defined via a system of generators.
Then $\Jg_0^\flat$ is generated by their images in $\Sg_0$, so it is even easier to control  $\Jg_0^\flat$
than $\Jg_0$.
\end{remark}
\section{Quantization of simple algebraic groups}
\label{secQSAG}
In the present section, we construct quantization of a special Poisson bracket
on a simple algebraic group. The quantized ring of polynomial functions is
realized as a quotient of the so called reflection equation algebra and simultaneously as a subalgebra
in the quantized universal enveloping algebra.
\subsection{Simple groups as Poisson Lie manifolds}
Let $\g$ be a complex simple Lie algebra and $G$ be a connected Lie group corresponding to
$\g$.
An element $g\in G$ is called semisimple if it belongs to a maximal torus in $G$.

Let $r$ be a classical $r$-matrix defining a factorizable Lie bialgebra structure on $\g$.
Denote by $r_-$ and $\Omega$ its skew and symmetric parts, respectively.
We assume $r$ to be normalized so that $\Omega$ is the inverse (canonical
element) of the Killing form on $\g$.

Given an element $\xi\in \g$ let $\xi^l$ and $\xi^r$ denote, respectively the left- and right invariant
vector fields on $G$ generated by $\xi$:
$$
(\xi^l f)(g)=\frac{d}{dt}f(ge^{t\xi})|_{t=0},
\quad
(\xi^r f)(g)=\frac{d}{dt}f(e^{t\xi}g)|_{t=0}
$$
for every smooth function $f$ on $G$.

The Semenov-Tyan-Shansky (STS) Poisson structure \cite{STS} on the group $G$ is defined
by the bivector field
\be
r_-^{l,l}+r_-^{r,r}-r_-^{r,l}-r_-^{l,r}
+\Omega^{l,l}-\Omega^{r,r}+\Omega^{r,l}-\Omega^{l,r}
=
r_-^{\ad,\ad} +(\Omega^{r,l}-\Omega^{l,r}).
\label{STSbr}
\ee
Here $\xi^{\ad}:=\xi^l-\xi^r$, $\xi\in \g$.

Consider the group $G$ as a $G$-space with respect to
the adjoint action.
The STS Poisson structure makes $G$ a Poisson-Lie manifold
over the Poisson Lie group $G$ endowed with the bracket $r^{l,l}-r^{r,r}$.
In fact, $G$ is a Poisson Lie manifold not only over $G$, but over
a Poisson Lie group corresponding to the double Lie bialgebra $\D\g$.
\begin{thm}[\cite{AM}]
Symplectic leaves of the  STS Poisson structure coincide with conjugacy classes in $G$.
\end{thm}

Let us compute the restriction of the STS bracket to the
class $C_g$ of a semisimple element $g\in G$.
The Lie algebra $\g$ splits into the direct sum
$\g=\l\oplus\m$ of vector spaces, where $\l$ is the eigenspace of $\Ad(g)$ corresponding to the eigenvalue
$1$ and $\m$ is the $\Ad(g)$-invariant subspace where $\Ad(g)-\id$ is invertible.
This decomposition is orthogonal with respect to the Killing form.
The tangent space to $C_g$ at the point $g$ is identified with $\m$,
while $\l$ is the Lie algebra of the centralizer of $g$.

Let $\{\xi_\mu\}\subset \m$ be an orthonormal base of eigenvectors of $\Ad(g)$
labelled by the eigenvalues of $\Ad(g)$. We have $(\xi_\mu,\xi_\nu)=0$
unless $\mu\nu=1$ and assume the normalization $(\xi_\mu,\xi_{\mu^{-1}})=1$.
One can check that the restriction of the STS bracket to the tangent spaces at the point $g$
is the bivector
$$r_{\m\wedge \m}+ \sum_{\mu}
\frac{\mu+1}{\mu-1}
\xi_\mu \tp
\xi_{\mu^{-1}}\in \m\wedge \m,
$$
where the first term is the projection of $r$ to $\m\wedge \m$.
The second term is correctly
defined since $\Ad(g)-\id$  is invertible on $\m$.
\subsection{Quantization of the STS bracket on the group}
\label{ssecqSTS}
In this subsection we describe quantization of the STS bracket
on classical matrix groups in terms of generators and relations.

By $\C_\hbar[G]$ we denote the $\U_\hbar(\D\g)$-equivariant deformation (quantization) of the affine
ring $\C[G]$. This quantization in the form of star product was constructed in \cite{DM3}.
Below we give a description of $\C_\hbar[G]$ in terms of generators and relations,
using the so called reflection equation algebra.

Let $(V,\pi)$ be the defining representation of $\U_\hbar(\g)$ and $R$ the
image of $\Ru$ in $\End(V^{\tp 2})$.
Denote by $\Kc$ the quotient of the tensor algebra of $\End(V^*)$ by the quadratic relations
\be
R_{21}K_1 R_{12} K_2 = K_2 R_{21} K_1 R_{12}.
\label{RE}
\ee
Here $K=||K^i_j||$ is the matrix of the generators forming the standard matrix base in
$\End(V^*)$. The algebra $\Kc$ is called reflection equation (RE) algebra, \cite{KSkl}.
It is a $\C[[\hbar]]$-flat deformation of the polynomial ring
$\C[\End(V)]$ only for $\g=sl(n)$.

The $\U_\hbar(\g)$-equivariant quantization of $\C[G]$ can be realized as a quotient
of $\Kc$. Below we describe the corresponding ideals only for $G$ orthogonal
and symplectic.
That can also be done for the  case  $G=SL(n)$, but then the algebra $\Kc$  is good enough for our
purposes, so it is even more convenient to work with $\Kc$ rather than $\C_\hbar[SL(n)]$.

Assume $\g$ to be orthogonal or
symplectic. Let $B$ be the $\U_\hbar(\g)$-invariant element from $V\tp V$, \cite{FRT}.
\begin{propn}
The algebra $\C_\hbar[G]$ is isomorphic to the quotient of $\Kc$ by the $\U_\hbar(\g)$-invariant ideal of  relations
\be
\label{orhtog_re}
R_1^tK^t\bigr ( (R^t_{1'})^{-1}B^{-1} ( R_{2'})^{-1}\bigl)R_2 K= B^{-1},
\quad
K R_1 B K^tR_2^t =R_{1'}B R_{2'}^t.
\ee
\end{propn}
A proof that this quotient is $\C[[\hbar]]$-free is given in \cite{M}.

\subsection{Embedding of $\C_\hbar[G]$ in $\U_\hbar(\g)$}

Recall that a (left) Yetter-Drinfeld (YD) module over $\U_\hbar(\g)$ is simultaneously a left
$\U_\hbar(\g)$-module, a left $\U_\hbar(\g)$-comodule, and these two structures are compatible
in a certain way, \cite{Y}. A YD algebra over  $\U_\hbar(\g)$ is a  $\U_\hbar(\g)$-module algebra,
$\U_\hbar(\g)$-comodule algebra and a YD module with respect to these structures.
An example of YD module algebra is $\U_\hbar(\g)$ considered
as the adjoint module and comodule via the coproduct.

Let us describe a method of producing YD modules (module algebras)
out of $\U_\hbar(\D\g)$-modules, where $\D\g$ is the double of $\g$.
For a factorizable semisimple Lie bialgebra $\g$ the algebra $\U_\hbar(\D\g)$
is isomorphic to the twisted tensor square \tw{\U_\hbar(\g)}{\Ru}{\U_\hbar(\g)}.
The universal R-matrix $\hat\Ru$ of $\U_\hbar(\D\g)$
is expressed through the universal R-matrix $\Ru$ of $\U_\hbar(\g)$ by the formula
$\hat\Ru=\Ru_{41}^{-1}\Ru_{31}^{-1}\Ru_{24}\Ru_{23}$. It is easy to see that the right
component of  $\hat \Ru$ lies, in fact, in $\U_\hbar(\g)$ (via diagonal embedding).
Then any $\U_\hbar(\D\g)$-module (module algebra) becomes a YD module (YD module algebra)
when equipped with the $\U_\hbar(\g)$-coaction
$\delta(a) = \hat\Ru_2\tp \hat\Ru_1\tr a$. The $\U_\hbar(\g)$-action is induced through
the embedding of $\Delta\colon\U_\hbar(\g)\hookrightarrow \U_\hbar(\D\g)$.

The algebra $\C_\hbar[G]$ is a $\U_\hbar(\D\g)$-algebra and hence
YD module algebra, by the above construction.
It can be realized in
as a $\C[[\hbar]]$-submodule  in the FRT dual to $\U_\hbar(\g)$.
The Hopf pairing between the FRT dual and $\U_\hbar(\g)$ induces
a paring between $\C_\hbar[G]$ and  $\U_\hbar(\g)$. By means of this pairing,
the universal RE matrix $\Q$ implements a $\U_\hbar(\g)$-algebra homomorphism
\be
\label{REemb}
\C_\hbar[G]\to \U_\hbar(\g), \quad
a \mapsto \langle a, \Q_1\rangle \Q_2.
\ee
It is easy to check using the explicit form of the $\Ru$ matrix
that $\C_\hbar[G]$ lies in $\breve \U_\hbar(\g)\subset \U_\hbar(\g)$.
The map (\ref{REemb}) is, in fact, a homomorphism of  YD algebras.

\begin{propn}
\label{q-n-gr}
The map (\ref{REemb}) is  embedding.
\end{propn}
\begin{proof}
The proof easily follows from  Proposition \ref{Lemma_aux} after a slight
adaptation to YD module algebras.
The algebra $\C_\hbar[G]$ decomposes into the direct sum
$\oplus_{W_0} W\tp W^*$ taken over the simple finite dimensional $\g$-modules.
Each summand is a YD module and its quotient $\mod \hbar$ is $\D\g=\g\oplus \g$-irreducible
with multiplicity one.

We have a pair of YD module algebras $\Sg=\C_\hbar[G]$ and $\Tg=\U_\hbar(\g)$. The former
is admissible (as a YD-module), while the latter has no torsion.
Set $\Psi$ to be  the map (\ref{REemb}). Its image and kernel are YD modules and
are free over $\C[[\hbar]]$. Therefore $\Sg$ decomposes into the direct sum
$\ker \Psi\oplus \im \Psi$. The rest of the proof is readily adapted from the proof of Proposition \ref{Lemma_aux}
if one observes that YD modules become $\D\g$-modules in the quasi-classical limit
(the first order in $\hbar$). Put $\Jg=0$. The algebra $\Sg_0=\C[G]$ has no non-zero $\D\g$-invariant
ideals. Thus we conclude $\ker \Psi=\Jg =0$.
\end{proof}
\section{Center of the algebra $\C_\hbar[G]$}
\label{subsecMC}

Let $G$ be a simple complex algebraic group.
If $G$ is simply connected, then $\C[ G]$ is a free module over the subalgebra of invariants  $I( G)$,
\cite{R}.
More precisely, there exists a submodule $\E_0\subset \C[ G]$ such that the multiplication
map $I( G)\tp \E_0\to \C[ G]$ is an isomorphism of vector spaces.
Each isotypic component in $\E_0$ enters with  finite multiplicity. This fact has a quantum analog.

\begin{thm}[\cite{M}]
\label{G=IE}
Let $G$ be a simple complex algebraic group and let $\C_\hbar[ G]$ be the
$\U_\hbar(\D\g)$-equivariant quantization of $\C[G]$ along
the STS bracket. Then
\\
i) the subalgebra $I_\hbar(G)$ of $\U_\hbar(\g)$-invariants coincides with the center of $\C_\hbar[ G]$,
\\
ii) $I_\hbar(G)\simeq I(G)\tp \C[[\hbar]]$ as $\C$-algebras.\\
Let $\hat G$ be the simply connected covering of $G$.
Then
\\
iii) $\C_\hbar[\hat  G]$ is a free $I_\hbar(\hat G)$-module
generated by a $\U_\hbar(\g)$-submodule $\E\subset \C_\hbar[\hat  G]$. Each isotypic component in
$\E$ is $\C[[\hbar]]$-finite.
\end{thm}
Theorem \ref{G=IE} i) implies, in particular,
that the center of $\C_\hbar[G]$ is the intersection of $\C_\hbar[G]$  with the center
of $\U_\hbar(\g)$.

If $G$ is not simply connected, the Theorem \ref{G=IE} iii) will be true only if the classical
algebra of invariants $I(G)$ is polynomial. That is the case, e.g. for $G=SO(2n+1)$,
however not so for $G=SO(2n)$.

We will use Theorem \ref{G=IE} for quantization of conjugacy classes  in Subsection \ref{ssecQThnm}.
Important for us will be the following fact.
\begin{propn}
\label{quotient}
Let $G$ be a simple complex algebraic group.
Suppose $\la$ is a character of $I_\hbar(G)$, i. e. a unital homomorphism to $\C[[\hbar]]$.
Denote by $J_\la$ the ideal in $\C_\hbar[G]$ generated by $\ker\la$.
Then the quotient $\C_\hbar[G]/J_\la$ is an admissible $\U_\hbar(\g)$-module.
\end{propn}
\begin{proof}
The group $G$ is a quotient of its universal covering group $\hat G$ over
a finite central subgroup $Z\subset \hat G$. The affine ring $\C[G]$ is embedded
in $\C[\hat G]$ as a subalgebra of $Z$-invariants with respect to
the regular action. The group $Z$ naturally acts on the quantized
algebra $\C_\hbar[\hat G]$, and the subalgebra of $Z$-invariants
is exactly $\C_\hbar[G]$. Accordingly, $I_\hbar(G)$ is the subalgebra
of $Z$-invariants in $I_\hbar(\hat G)$. The latter is finitely generated
over $I_\hbar(G)$, since $Z$ is finite.

Let $\la$ be a character of $I_\hbar(G)$. Denote by $J_\la$ and $\hat J_\la$
the ideals  generated by $\ker \la$ in  $\C_\hbar[G]$ and $\C_\hbar[\hat G]$,
respectively. We have $J_\la=\hat J_\la\cap \C_\hbar[G]$, as follows from $Z$-invariance.
Therefore the quotient  $\C_\hbar[G]/J_\la$ is embedded in
$\C_\hbar[\hat G]/\hat J_\la$. The quotient $\hat I_\hbar(G)/\hat J_\la$ is
$\C[[\hbar]]$-finite and therefore the module $\C_\hbar[\hat G]/\hat J_\la$ has finite isotypic components.
Its submodule $\C_\hbar[G]/J_\la$ has also finite components. It is admissible, being a submodule
of an admissible module.
\end{proof}

Let us emphasize that the quotient $\C_\hbar[G]/J_\la$ is $\C[[\hbar]]$-free if the classical subalgebra of invariants
is a polynomial algebra, cf. remark after Theorem \ref{G=IE}. Otherwise it may have $\hbar$-torsion.
\subsection{The case of classical series}
\label{theta}
The case of $G=SO(2n)$ differs from other classical matrix groups. Let us focus
on a more simple case of the $A,B,C$ series first.
\begin{propn}
\label{center}
The elements $\tau^\ell$, $\ell=1,\ldots ,N$, generate the $\hbar$-adic completion of the center in
the algebras $\C_\hbar[SL(n)]$, $\C_\hbar[SO(2n+1)]$, and $\C_\hbar[Sp(n)]$.
\end{propn}
\begin{proof}
The center in $\C_\hbar[G]$ is isomorphic to classical subalgebra of invariants extended by $\C[[\hbar]]$,
as stated in  Theorem \ref{G=IE}.
Modulo $\hbar$ it is generated by the classical limits of $\tau^\ell$, i.e.
the traces of matrix powers. This implies the statement.
\end{proof}

In the rest of the subsection we assume $G=SO(2n)$.
Let $V$ be the defining $\U_\hbar\bigl(so(2n)\bigr)$-module
and $\{\pm \ve_i\}_{i=1}^n=\Lambda(V)$ be the set of weights.

In the classical limit, the  $n$-th exterior power of $V_0$ splits into two irreducible submodules
with highest weights $\pm\ve_n+\sum_{i=1}^{n-1} \ve_i $, \cite{VO}.
Therefore the $\U_\hbar\bigl(so(2n)\bigr)$-module $\wedge^n_q V$ (q-anti-symmetrized)
is a direct sum of two modules $(W_\pm,\pi_\pm)$ with the highest weights
$\pm\ve_n+\sum_{i=1}^{n-1} \ve_i $, respectively.
Let $p_\pm\colon V^{\tp n}\to W_\pm$ be the intertwining projectors.
The element $(\Delta^{n}\tp \id)(\Q)$ is expressed through $\Q$ and the
numerical matrix $R=(\pi\tp \pi)(\Ru)$ in the defining representation, by virtue of
(\ref{dQ1}).
Hence  $\Q_{W_\pm}$ can be  explicitly expressed through the matrix $\Q_{V}$ whose entries generate $\C_\hbar[G]$
within $\U_\hbar(\g)$.
In the same fashion, we can define the invariant matrices $K_\pm\in \End(W_\pm)\tp \Kc$
expressing them through $K$, $R$, and $p_\pm$ by the same formulas.
The matrices $\Q_{W_\pm}$ are obtained from $K_\pm$ via the projection $\Kc\to \C_\hbar[G]$.

Define a central element $\tau^-$  of $ \U_\hbar\bigl(so(2n)\bigr)$ by setting
\be
\tau^-:=\Tr_q(\Q_{W_+})-\Tr_q(\Q_{W_-}).
\ee
Similarly we introduce the central elements $\Tr_q(K_+)-\Tr_q(K_-)\in \Kc$.
\begin{propn}
Let $\la\in \h^*$ be a weight. Then
$
\chi^\la(\tau^-)=\prod_{i=1}^n(q^{2(\la+\rho,\ve_i)}-q^{-2(\la+\rho,\ve_i)}).
$
\end{propn}
\begin{proof}
\newcommand{\Ch}{\mathrm{Ch}}
By Corollary \ref{Tr-Tr}, the central element $\Tr_q(\Q_{W_\pm})$ acts on a module with highest weight
$\la$ as multiplication by the scalar $\Tr\bigl(\pi_\pm(q^{2(h_\la+h_\rho)})\bigr)$. This scalar becomes
the group character $\ch_\pm(t)$ (a class function) associated with $W_\pm$ upon the substitution $q^{2(\la+\rho,\ve_i)}\to t_i$,
where $\{t_i\}_{i=1}^n$ are the coordinate functions on the maximal torus in $SO(2n)$.

It is known \cite{We} that the ring of characters of the group $SO(2n)$ is isomorphic
to $\C[t_1,\ldots,t_n,t_1^{-1},\ldots,t_n^{-1}]^{\Wm_{SO(2n)}}$.
The group $\Wm_{SO(2n)}$ acts by permutations of the pairs $(t_i,t^{-1}_i)$
and even number of inversions $t_i\leftrightarrow t_i^{-1}$.

The difference $\ch_+(t)-\ch_-(t)$ changes sign under every inversion $t_i\leftrightarrow t_i^{-1}$,
hence it is divided by $\prod_{i=1}^n (t_i-t_i^{-1})$. In fact, $\ch_+(t)-\ch_-(t)$
equals $\prod_{i=1}^n (t_i-t_i^{-1})$.
This is verified in a standard way by comparing
the highest and lowest terms with respect to a natural lexicographic
ordering in $\C[t_1,\ldots,t_n,t_1^{-1},\ldots,t_n^{-1}]$
(and elementary analysis of the weight structure of $\Lambda^n V_0$).
\end{proof}
As a corollary, we obtain the following.
\begin{propn}
The elements $\tau^-$ and $\tau^\ell$, $\ell=1,\ldots ,N$, generate the $\hbar$-adic completion of the
center of  $\C_\hbar[SO(2n)]$.
\end{propn}
\begin{proof}
Similarly to Proposition \ref{center}.
\end{proof}

\section{Quantization of conjugacy classes}
\label{secQCC}
\subsection{Non-exceptional classes}
\label{secSSCC}
Semisimple conjugacy  classes are the only closed classes of a simple
algebraic group, \cite{S}.
Among semisimple classes of a classical matrix group $G$ we select those which
are isomorphic to semisimple orbits in $\g^*$ as affine algebraic varieties; we call
them "non-exceptional".
The isomorphism is implemented by the Cayley transformation $X\mapsto (1 \mp X)(1 \pm X)^{-1}$ for $G$ orthogonal and symplectic.
Thus all such classes have Levi subgroups as stabilizers.
Non-exceptional covers all semisimple classes
for $SL(n)$ and all classes with Levi stabilizers for $Sp(n)$. For  orthogonal groups
non-exceptional are classes of matrices with no eigenvalues $+1$
and $-1$ simultaneously.

Let $g\in G$ be a semisimple element. It satisfies a matrix polynomial equation
$p(g)=0$ and defines a character $\chi^g$ of the subalgebra of invariants in $\C[G]$.
\begin{thm}
Defining ideal $\Nc(C_g)$ of a non-exceptional conjugacy class $C_g$ is generated by
the kernel of $\chi^g$ and by the entries of
the minimal matrix polynomial for $g$.
\end{thm}
A proof of this theorem will be given elsewhere. It is based on the following facts.
a) Non-exceptional classes are isomorphic
to semisimple coadjoint orbits via the Cayley transformation. b) The defining ideals of semisimple
coadjoint orbits can be obtained
as classical limits of annihilators of generalized Verma modules, \cite{DGS}.
c) It was shown in \cite{G} (see also \cite{J2}) that for certain weights the annihilator of the generalized Verma module
is generated by a copy of adjoint module in $\U(\g)$ and the kernel of the central character.
This is sufficient to describe the defining  ideals for all semisimple coadjoint orbits and
hence for non-exceptional conjugacy classes.

\subsection{The quantization theorem}
\label{ssecQThnm}
In this subsection $G$ is a simple complex algebraic group from the classical series and $\g$ its Lie algebra.
Fix a non-exceptional semisimple element $g$ in the maximal torus $T\subset G$.
Let $\l$ be the Lie algebra of the centralizer of $g$, which is a Levi subgroup.
Take $h_\la\in \h$ such that $e^{h_\la}=g$ and $\la\in \c^*_{reg}$ is a regular
character of $\l$.

As was stated in Remark \ref{breve}, the generalized Verma modules can be
naturally considered over the ring of polynomials $\C[\c^*]$, where $\c$ is the
center of $\l$. In this case, the parabolic induction is performed
from the representation of $\breve\U_\hbar(\l)$ in $\C[\c^*][[\hbar]]$.
Propositions  \ref{spec} and \ref{Spec_parabolic} are valid for such modules, since they are
valid for generic element of $\c^*$.
Let $(c_i)$ be coordinate functions generating $\C[\c^*]$
and $(\la_i)$  the coordinates of a regular element $\la\in \c^*_{reg}$.
The homomorphism map ${c_i}\mapsto \la_i/(2\hbar)$ defines a character of
$\breve \U_\hbar(\l)$ (but not of $\U_\hbar(\l)$). Let $\breve M_{\la/2\hbar}$
denote the corresponding Verma module over $\breve \U_\hbar(\g)$
 It is easy to see that the $\breve \U_\hbar(\g)$-action
on $\End(\breve M_{\p,\la/2\hbar})$ extends to an action of $\U_\hbar(\g)$.

Recall that $(V,\pi)$ stands for the defining representation of $\U_\hbar(\g)$.
Also recall from Subsection \ref{ssecqSTS} that the entries of the matrix
$\Q_V=\pi(\Q_1)\tp \Q_2\in \End(V)\tp \U_\hbar(\g)$
generate the subalgebra $\C_\hbar[G]\subset \U_\hbar(\g)$, which is the
equivariant quantization of $\C[G]$.
\begin{thm}
\label{main}
For non-exceptional $\la\in \c^*\subset \h^*$ the image of
$\C_\hbar[G]$ in $\End(\breve M_{\p,\la/2\hbar})$  is a quantization of the ring of regular
functions on the conjugacy class $C_{\exp(h_\la)}\subset G$.
The ideal of the quantized class is generated by the entries of the minimal polynomial
in $\Q_{V}$ over the kernel of a character of $I_\hbar(G)$.
\end{thm}
\begin{proof}
The composite map
\be
\label{quantizing_seq}
\C_\hbar[G]\hookrightarrow \breve\U_\hbar(\g) \to \End(\breve M_{\p,\la/2\hbar})
\ee
is a $\U_\hbar(\g)$-equivariant algebra homomorphism.
When restricted to the center in $\C_\hbar[G]$, this map defines a character
of $I_\hbar(G)$.
Indeed, a central element from $\C_\hbar[G]$ acts on
$\breve M_{\p,\la/2\hbar}$ as multiplication by a scalar. That scalar
is a polynomial in the eigenvalues of the matrix $\Q_V$, as follows from Proposition \ref{trace_auxi}.
The eigenvalues are given by (\ref{spec_V}), where $\la$ should be replaced
by $\la/2\hbar$. This proves that the representation in $\breve M_{\p,\la/2\hbar}$
defines a character of the center of $\breve \U_\hbar(\g)$ and therefore
a character of $I_\hbar(G)$. This character is a deformation of the $I(G)$-character $\chi^g$,
and will be denoted by $\chi^g_\hbar$.

The matrix $\Q_V$ satisfies a polynomial equation for a polynomial
$p$ with simple roots. The roots of $p$ are given by Theorem \ref{Spec_parabolic}, where again
$\la$ should be replaced by $\la/2\hbar$. These roots are regular in $\hbar$ and go
over to the eigenvalues $e^{(\la,\ve_i)}$ of the matrix $g=e^{h_\la}\in G$.

Next we are going to apply Proposition \ref{Lemma_aux}.

Put $\Sg$ to be the quotient of $\C_\hbar[G]$ over the ideal generated by $\ker \chi^g_\hbar$.
The algebra $\Sg$ is admissible, due to Proposition \ref{quotient}.
Put $\Tg=\End(\breve M_{\p,\la/2\hbar})$. The $\C[[\hbar]]$-module $\Tg$ has no torsion,
as a subspace of endomorphisms of a torsion free module.
By construction, the map (\ref{quantizing_seq}) factors through
an equivariant algebra map $\Psi \colon \Sg\to \Tg$.

Put $\Jg \subset \C_\hbar[G]$ to be the ideal generated by
the entries of $p(\Q_V)$ projected to $\Sg$.
By construction, $\Jg\subset \ker \Psi$. The image $\Jg_0^\flat$ of $\Jg_0$ in
$\Sg_0=\C[G]/\bigl(I(G)-\chi^g\bigr)$ is a maximal $\g$-invariant ideal, because it is the image
of the maximal $\g$-invariant ideal $\Nc(C_g)\subset \C[G]$.

Thus the conditions of Proposition \ref{Lemma_aux} are satisfied.
Therefore  $\Sg/\Jg$ is the quantization of $\C[C_g]$ and $\ker \Psi=\Jg$.
In terms of $\C_\hbar[G]$, the defining ideal of $\C[C_g]$ is generated by
the kernel of $\chi^g_\hbar$ and the entries of $p(\Q_V)$.
\end{proof}
Next we compute the quantized ideals of non-exception conjugacy classes.

\subsection{Ideals of quantized non-exceptional conjugacy classes}
First of all, we specialize the formula (\ref{char_V}) for $\g$ being the simple matrix Lie
algebra and $V$ the defining representation  of $\U_\hbar (\g)$.
Since every weight has multiplicity one, we rewrite (\ref{char_V}) as
\be
\label{char_V1}
\chi^\la(\tau^\ell)&=&
\sum_{\nu_i\in \Lambda(V)} x_{\nu_i}^\ell
\prod_{\al\in \Rm_+}\frac{q^{(\la+\rho+\nu_i,\al)}-q^{-(\la+\rho+\nu_i,\al)}}{q^{(\la+\rho,\al)}-q^{-(\la+\rho,\al)}}
.
\ee

Here $x_{\nu_i}$, $\nu_i\in \Lambda(V)$, are the roots of the characteristic polynomial
for $\Q$ considered as an operator on $V\tp M_{\la}$.
They are related to the highest weight $\la$ by $x_{\nu_i}=q^{2(\la+\rho,\nu_i)-2(\rho,\nu)}$, $\nu_i\not=0$
and $x_0=q^{-2(\rho,\nu)-(\nu,\nu)}$.

Perform the following substitution in (\ref{char_V1}).
Set $x_0:= q^{-2n}$ for $\g=so(2n+1)$, $x_i=x_{\nu_i}$ for all $\g$, and
$x_{i'}=x_{-\nu_i}=x^{-1}_i q^{-4(\rho,\nu)}$
for $\g$ orthogonal and symplectic. Here $i$ ranges from $1$ to $n=\rk\; \g$.
The eigenvalue $x_0$ is present
only for $\g=so(2n+1)$.
Below we use the convention $i'=N+1-i$, cf. Subsection \ref{ssecVR}.
Under the adopted enumeration of weights of $V$ the
eigenvalues $x_i$ and $x_{i'}$, $i=1,\ldots, N$, correspond to the weights of the opposite sings.
 As a result of this substitution, we obtain the functions
\be
\vt^{(\ell)}_{sl(n)}(x)
&=&\sum_{i=1}^n x_i^\ell \prod_{j=1 \atop j\not =i}^n\frac{q x_i-x_j\bar q}{x_i- x_j},
\nn
\\
\vt^{(\ell)}_{so(2n+1)}(x)
&=&
\sum_{i=1}^n x_i^\ell \;\frac{q x_i-x_0 }{x_i - x_0q}
\prod_{j=1 \atop j\not =i}^n
\frac{q x_i-x_j\bar q}{x_i- x_j}\;
\frac{q x_i-x_{j'}\bar q}{x_i- x_{j'}}
\>+\>i\leftrightarrow i'
+x_0^{\ell}
,
\nn
\\
\vt^{(\ell)}_{sp(n)}(x)
&=&\sum_{i=1}^n x_i^\ell\; \frac{q^2 x_i-x_{i'}\bar q^2}{x_i- x_{i'}}
\prod_{j=1 \atop j\not =i}^n
\frac{q x_i-x_j\bar q}{x_i- x_j}\;
\frac{q x_i-x_{j'}\bar q}{x_i- x_{j'}}
\> + \> i\leftrightarrow i'
,\nn
\\
\vt^{(\ell)}_{so(2n)}(x)
&=&\sum_{i=1}^n x_i^\ell\;\prod_{j=1 \atop j\not =i}^n\frac{q x_i-x_j\bar q}{x_i- x_j}\;
\frac{q x_i-x_{j'}\bar q}{x_i- x_{j'}}
\>+\>
i\leftrightarrow i',
\nn
\\
\vt^{-}_{so(2n)}(x)
&=&
q^{2n(\rho,\nu)}\prod_{i=1}^n\bigr(x_i -x_{i'} \bigl),
\nn
\ee
where $\bar q:=q^{-1}$.
Remark that $\vt^{(\ell)}_\g$ are in fact polynomial in $x_i, x_{i'}, x_0$.

Fix a pair $(\mub,\nb)\in \C^k\times \Z^k_+$ such that $\sum_{i=1}^k n_i=n$.
Define the vector
\be
x(\mub):=\bigl(\underbrace{\mu_1,\ldots, \mu_1}_{n_1},\ldots , \underbrace{\mu_k,\ldots, \mu_k}_{n_k}\bigr)
\in \C^n
\ee
and $x_q(\mub)$ obtained from $x(\mub)$ through replacing the
"constant" string $(\mu_i,\ldots,\mu_i)$
by "quantum" $(\mu_i,\mu_iq^{-2},\ldots,\mu_iq^{-2(n_i-1)})$
for each $i=1,\ldots,k$ .
Define functions
$\vt^{(\ell)}_{\g}(\nb,\mub)$ via the
substitution  $x=x_q(\mub)$, $x_0=q^{-2n}$, to $\vt^{(\ell)}_\g(x)$.

Recall from Subsection \ref{ssecVR} that weights of the defining representation can be expressed in
terms of the orthogonal set $\{\ve_i\}$. Next we specialize different types of Levi subalgebras in $\g$.

It is convenient to consider $gl(n)$ instead of $sl(n)$.
The Levi subalgebras in $gl(n)$ have the form $\l=\oplus_{i=1}^k gl(n_i)$.
The  irreducible $\l$-submodules in the defining $\g$-module $V_0$
are labelled with highest weights
\be
\bigl\{\ve_1,\ve_{n_1+1},\ldots, \ve_{n_1+\ldots + n_{k-1}+1}\bigr\},
\label{l_sl}
\ee

The Levi subalgebras in the orthogonal and symplectic algebras
can be represented as $\l=\oplus_{i=1}^{k}\l_i$, where
$\l_i=gl(n_i)$ for $i=1,\ldots,k-1$ and $\l_k$ equals
either $gl(n_k)$ or $so(2n_k+1)$, $sp(2n_k)$, $so(2n_k)$
for $\g$ being $so(2n+1)$, $sp(n)$, $so(2n)$, respectively.
That is, the latter option corresponds to $\l_k$ of the same type as $\g$.

In the case $\l_k=gl(n_k)$ the irreducible $\l$-submodules in  $V_0$ are
labelled with the highest weights
\be
\bigl\{\ve_1,\ve_{n_1+1},\ldots ,\ve_{n_1+\ldots + n_{k-1}+1};-\ve_{n_1+\ldots + n_{k}},\ldots ,-\ve_{n_1+n_2},-\ve_{n_1}\bigr\},
\ee
and the zero weight for $so(2n+1)$. Each simple $\l$-modules enters with its dual.

For  $\l_k$ of the same type as $\g$  the irreducible $\l$-modules in the vector representation $V_0$ are
labelled with the highest weights
\be
\bigl\{\ve_1,\ve_{n_1+1},\ldots ,\ve_{n_1+\ldots + n_{k-2}+1};
\ve_{n-n_k+1}; -\ve_{n_1+\ldots + n_{k-1}},\ldots ,-\ve_{n_1+n_2},-\ve_{n_1} \bigr\},
\ee

We assume that the components $\mu_i$ of $\mub$ are non zero and pairwise distinct. For $\g$ orthogonal and symplectic
we also assume $\mu_i\not =\mu_j^{-1}$ for $i\not = j$ and $\mu_i^2\not =1$ unless otherwise
stated.

Below $Q=\Q_{V}$ is the matrix whose entries generate $\C_\hbar[G]\subset \U_\hbar(\g)$.

\subsubsection{The $\U_\hbar\bigr(gl(n)\bigl)$-case}
\label{gl}
Put $g=\sum_{i=1}^ng_ie_{ii} \in GL(n)$ to be the diagonal matrix with
$g_i=x(\mub)_i$.
The quantized ideal $\Nc_\hbar(C_g)\subset \C_\hbar[SL(n)]$ of the class $C_g$ is generated by the relations
$$
\prod_{i=1}^k(Q-\mu_i)=0,
\quad
\Tr_q(Q^\ell)=\vt^{(\ell)}_{sl(n)}(\nb,\mub).
$$
This reproduces the results of \cite{DM2}, where the functions $\vt^{(\ell)}_{sl(n)}(\nb,\mub)$
are written out in a manifestly polynomial form. In the present form, the central characters of
the RE algebras are
calculated for the general case of Hecke symmetries in \cite{GS}.
\subsubsection{The $\U_\hbar\bigr(so(2n+1)\bigl)$-case}
\underline{1. The case $\l_k=gl(n_k)$}.
Set
$
g:=
g=\sum_{i=1}^{2n+1}g_ie_{ii}\in SO(2n+1)$
with $g_i=x(\mub)_i$ for $i=1,\ldots, n$, $g_{n+1}=1$, and $g_{i}=x(\mub)_{i'}^{-1}$ for $i=n+2,\ldots, 2n+1$.
The quantized ideal $\Nc_\hbar(C_g)\subset \C_\hbar[SO(2n+1)]$ of the class $C_g$ is generated by the relations
\be
(Q-q^{-2n})\prod_{i=1}^k(Q-\mu_i)(Q-\mu_i^{-1}q^{-4n+2n_i})=0,
&\hspace{-15pt}&
\Tr_q(Q^\ell)=\vt^{(\ell)}_{so(2n+1)}(\nb,\mub).
\nn
\ee
\underline{2.
The case $\l_k=so(2n_k+1)$}.
Set
$g$ as in the case $\l_k=gl(n_k)$ with $\mu_k=1$.
The quantized ideal $\Nc_\hbar(C_g)\subset \C_\hbar[SO(2n+1)]$ of the class $C_g$ is generated by the relations
\be
(Q-\mu_k)\prod_{i=1}^{k-1}(Q-\mu_i)(Q-\mu_i^{-1}q^{-n+2n_i})=0,
&\hspace{-15pt}&
\Tr_q(Q^\ell)=\vt^{(\ell)}_{so(2n+1)}(\nb,\mub),
\nn
\ee
where $\mu_k=q^{2(n_k-n)}$.
\subsubsection{The $\U_\hbar\bigr(sp(n)\bigl)$-case}
\underline{1. The case $\l_k=gl(n_k)$}.
Set
$
g=\sum_{i=1}^{2n}g_ie_{ii}\in Sp(n)$
with $g_i=x(\mub)_i$ for $i=1,\ldots, n$,  and $g_{i}=x(\mub)_{i'}^{-1}$ for $i=n+1,\ldots, 2n$.
The quantized ideal $\Nc_\hbar(C_g)\subset \C_\hbar[Sp(n)]$ of the class $C_g$ is generated by the relations
\be
\prod_{i=1}^k(Q-\mu_i)(Q-\mu_i^{-1}q^{-4n+2(n_i-1)})=0,
&\hspace{-15pt}&
\Tr_q(Q^\ell)=\vt^{(\ell)}_{sp(n)}(\nb,\mub).
\nn
\ee
\underline{ 2. The case $\l_k=sp(2n_k)$}.
Set
$g$ as in the case $\l_k=gl(n_k)$ with  $\mu_k=1$.
The quantized ideal $\Nc_\hbar(C_g)\subset \C_\hbar[Sp(n)]$ of the class $C_g$ is generated by the relations
\be
(Q-\mu_k)\prod_{i=1}^{k-1}(Q-\mu_i)(Q-\mu_i^{-1}q^{-4n+2(n_i-1)})=0,
&\hspace{-15pt}&
\Tr_q(Q^\ell)=\vt^{(\ell)}_{sp(n)}(\nb,\mub),
\nn
\ee
where $\mu_k=q^{2(n_k-n)}$.
\subsubsection{The $\U_\hbar\bigr(so(2n)\bigl)$-case}
\underline{1. The case $\l_k=gl(n_k)$}.
Set
$
g=\sum_{i=1}^{2n}g_ie_{ii}\in SO(2n)$
with $g_i=x(\mub)_i$ for $i=1,\ldots, n$,  and $g_{i}=x(\mub)_{i'}^{-1}$ for $i=n+1,\ldots, 2n$.
The quantized ideal $\Nc_\hbar(C_g)\subset \C_\hbar[SO(2n)]$ of the class $C_g$ is generated by the relations
\be
&\prod_{i=1}^k(Q-\mu_i)(Q-\mu_i^{-1}q^{-4n+2(n_i+1)})=0,
\\[6pt]
&\Tr_q(Q_{W_+})-\Tr_q(Q_{W_-})=\vt^-(\nb,\mub),
\quad
\Tr_q(Q^\ell)=\vt^{(\ell)}_{so(2n)}(\nb,\mub).
\nn
\ee
Without specializing the value of $\Tr_q(Q_{W_+})-\Tr_q(Q_{W_-})$ we get a quantization of
the intersection of a $SL(2n)$-class with the group $SO(2n)$. That intersection
is an $O(2n)$-class consisting
of two isomorphic (as Poisson Lie manifolds) $SO(2n)$-classes.

\underline{ 2. The case $\l_k=so(2n_k)$}.
Set
$g$ as in the case $\l_k=gl(n_k)$ with  $\mu_k=1$.
The quantized ideal $\Nc_\hbar(C_g)\subset \C_\hbar[SO(2n)]$ of the class $C_g$ is generated by the relations
\be
(Q-\mu_k)\prod_{i=1}^{k-1}(Q-\mu_i)(Q-\mu_i^{-1}q^{-4n+2(n_i+1)})=0,
\quad
\Tr_q(Q^\ell)=\vt^{(\ell)}_{so(2n)}(\nb,\mub),
\nn
\ee
where $\mu_k=q^{2(n_k-n)}$.

\begin{remark}
To express the quantized ideals in terms of the generators $\{K^i_j\}\subset \Kc$
one should replace $Q$ by $K$ in the formulas above and impose additional relations
(\ref{orhtog_re}) in case $\g$ is orthogonal or symplectic. If $\g=sl(n)$, one
can consider relations of Subsection \ref{gl} as those in $\Kc$. There are no
additional relations needed, and the eigenvalues may take arbitrary pairwise
distinct values. This case has been studied in a two parameter setting in \cite{DM2}.
\end{remark}

\appendix{}
\section{More on central characters}
In this appendix we derive some formulas for central characters of Drinfeld-Jimbo quantum groups.

The coproduct on the Cartan-Weyl generators reads
\be
\Delta(e_\al)=e_{\al}\tp 1 + q^{h_{\al}}\tp e_{\al}+\ldots,
\quad
\Delta(e_{-\al})=e_{-\al}\tp q^{-h_{\al}} + 1 \tp e_{-\al}+\ldots.
\ee
The omitted terms have lower root vectors in each tensor component.
\be
\Ru=q^{\Omega_\h} \Ru'(e_{-\al}\tp e_{\al}) = q^{\Omega_\h} \hspace{-7pt}\mod \U'_\hbar(\n^-)\tp \U'_\hbar(\n^+),
\label{Rmat1}
\ee
where $\Ru'(e_{-\al}\tp e_{\al})$ is a series in $e_{-\al}\tp e_{\al}$, $\al\in \Rm_+$ and
$\Omega_\h\in \h\tp \h$ is the inverse to the Killing form (the canonical element)
restricted to $\h$.

Let $W$ be a finite dimensional $\U_\hbar (\g)$-module. Let $X\in \End(W)\tp \U_\hbar (\g)$
be an invariant matrix, i. e. commuting with
$\pi(x^{(1)})\tp x^{(2)}$ for all $x\in \U_\hbar (\g)$.
It is known that the q-trace defined by
$$\Tr_{q}(X^n):=\Tr(\pi_W(q^{2h_\rho})X^n)$$
belongs to the center of $\U_\hbar (\g)$.
Below we will derive formula for $\chi^\la\bigl(\Tr_q(\Q^\ell_W)\bigr)$,
where $\chi^\la$ is a central character $\U_\hbar (\g)$ corresponding to $\la$.

For every weight $\la\in \h^*$ define a linear endomorphism $\theta_\la$
of $\U_\hbar(\g)$ setting
$$\theta_\la(x):= q^{2h_\la+2h_\rho-2(\rho,\nu)}\Ru_1x \Ru_2.$$
This endomorphism restricts to the subspace  of $\h$-invariants in $\U_\hbar(\g)$, since $\Ru$ is of zero weight.
\begin{propn}
\label{trace_auxi}
Let $(W,\pi_W)$ be a finite dimensional $\U_\hbar(\g)$-module with the highest weight $\nu$.
Then
\be
\chi^\la\bigl(\Tr_q(\Q^\ell_W)\bigr)
=
q^{-\ell(\nu,\nu)}
(\Tr_W\circ\pi_W)\bigl(\theta_\la^{\ell}(q^{2h_\rho})\bigr)v_\la.
\ee
\end{propn}
\begin{proof}
We have
$$
\chi^\la\bigl(\Tr_q(\Q^\ell_W)\bigr)=\Tr_q\Bigl((\pi\tp\id)\bigl(\Delta(\zz^\ell)\bigr)\Bigr)
q^{-\ell\bigl((\la,\la)+2(\rho,\la)+(\nu,\nu)+2(\rho,\nu)\bigr)}.
$$
Using induction in $\ell$, we find
$$
\zz^\ell=(q^{2h_\rho}\upsilon)^\ell=q^{2\ell \:h_\rho}\upsilon^\ell=
q^{2\ell \:h_\rho}
\underbrace{\gm^{2\ell}(\Ru_{1'})\Bigl(\ldots \bigl(\gm^{2}(\Ru_1)}_\ell
\underbrace{\Ru_2\bigr)\ldots\Bigr)\Ru_{2'}}_\ell.
$$
Let $v_\la$ be the highest weight vector of $M_\la$.
Consider a  linear map $\psi_\la\colon \U_\hbar(\g)\to \C[[\hbar]]$
defined by $x v_\la=\psi_\la (x)v_\la + \mbox{lower weight terms}$.
Then
$$
\chi^\la\bigl(\Q^\ell_W\bigr)= q^{-\ell(\nu,\nu)}(\Tr\circ\pi)(\Upsilon).
$$
where we set $\Upsilon=q^{-\ell(\la,\la)-2\ell(\rho,\la)-2\ell(\rho,\nu)}q^{2h_\rho}(\id \tp  \psi_\la)\bigl(\Delta(\zz^\ell)\bigr)\in \U_\hbar(\g)$.
Using the expression (\ref{Rmat1}) for the universal R-matrix and formulas for comultiplication on generators,
we compute $\Upsilon$ to be
$$
q^{-2\ell(\rho,\nu)}q^{2h_\rho}q^{2\ell \:h_\rho}
\underbrace{\{q^{h_\la}q^{h_\la}\gm^{2\ell}(\Ru_{1'})q^{-h_\la}\}\ldots
\{q^{h_\la}q^{h_\la}\gm^2(\Ru_1)q^{-h_\la }\}}_\ell
\underbrace{\{q^{h_\la}\Ru_2\}\ldots\{q^{h_\la}\Ru_{2'}\}}_\ell
.
$$
We can push the factor $q^{h_\la}$ from before each copy $\Ru_2$ to the left till it meets the factor
$q^{-h_\la}$ following the corresponding component $\gm^{2j}(\Ru_1)$. Since $\Ru$ is of zero weight,
$q^{h_\la}$ commutes with the factors in between.
Now recall that $\gm^2$ is given by conjugation
with $q^{-2h_\rho}$. Taking this into account, we obtain
\be
\Upsilon=
q^{-2\ell(\rho,\nu)}
\underbrace{\{q^{2h_\la+2h_\rho}\Ru_{1'}\}\ldots
\{q^{2h_\la+2 h_\rho}\Ru_1\}}_\ell
\;q^{2 h_\rho}
\underbrace{\{\Ru_2\}\ldots\{\Ru_{2'}\}}_\ell
=
\theta^{\circ\ell}(q^{2 h_\rho})
.
\nn
\ee
This makes the proof immediate.
\end{proof}
\begin{corollary}
\label{Tr-Tr}
Let $W$ be a finite dimensional $\U_\hbar(\g)$-module.
Then
$$
\chi^\la\bigl(\Tr_q(\Q_{W}) \bigr)=
\Tr\bigl(\pi_W(q^{2(h_\la+h_\rho)})\bigr).
$$
\end{corollary}
\begin{proof}
Notice that $\Ru_1q^{2h_\rho}\Ru_2$ equals $z$.
\end{proof}

\end{document}